# Some nice Hankel determinants

J. CIGLER

Fakultät für Mathematik, Universität Wien

*johann.cigler@univie.ac.at*

http://homepage.univie.ac.at/johann.cigler

**Abstract**

I study Hankel determinants of a class of sequences which can be interpreted as generalizations of the Catalan numbers and the central binomial coefficients. They follow a modular pattern with a frequent appearance of zeroes, so that the theory of orthogonal polynomials is not applicable. Even so our theorems and conjectures show some similarity with the relations between Catalan numbers and Fibonacci polynomials or between central binomial coefficients and Lucas polynomials.

**1. A survey of some known results**

As is well known there is a close relationship between orthogonal polynomials and Hankel determinants.

The **Hankel determinants** of order $k$ of a sequence $(a_n)_{n \geq 0}$ are the determinants

$$d_k(n) = \det\left(a_{i+j+k}\right)_{i,j=0}^{n-1}. \quad (1.1)$$

Let $F$ be a **linear functional** on the vector space $\mathbb{C}[z]$ of polynomials in the variable $z$. A sequence $(p_n)_{n \geq 0}$ of monic polynomials $p_n(z)$ with $\deg p_n = n$ is called **orthogonal with respect to** $F$ if $F(p_m p_n) = r_n [n = m]$ for some numbers $r_n \neq 0$. The numbers $\mu_n = F(z^n)$ are called **moments** of $F$.

The monic orthogonal polynomials $p_n(z)$ and their moments $\mu_n$ are related by the formula

$$p_n(z) = \frac{1}{d_0(n)} \det \begin{pmatrix} \mu_0 & \mu_1 & \cdots & \mu_{n-1} & 1 \\ \mu_1 & \mu_2 & \cdots & \mu_n & z \\ \mu_2 & \mu_3 & \cdots & \mu_{n+1} & z^2 \\ \vdots & & & & \vdots \\ \mu_n & \mu_{n+1} & \cdots & \mu_{2n-1} & z^n \end{pmatrix}. \tag{1.2}$$

For these polynomials are monic and $F(z^k p(z)) = 0$ for $0 \leq k < n$.

By Favard's formula for each sequence of monic orthogonal polynomials there exist numbers $s(n), t(n)$ such that $p_n(z)$ satisfies a 3-term recurrence

$$p_n(z) = (z - s(n-1))p_{n-1}(z) - t(n-2)p_{n-2}(z) \tag{1.3}$$

with initial values $p_{-1} = 0$ and $p_0 = 1$.

Then the Hankel determinants of order 0 and 1 of the moment sequence are $d_0(n) = \prod_{i=1}^{n-1} \prod_{j=0}^{i-1} t(j)$ and $d_1(n) = (-1)^n p_n(0) d_0(n)$.

In order to make this paper self-contained proofs will be given in section 2.

Most results of this paper can be interpreted as extensions of the following simple examples. (Precise definitions will be given later).

**Example 1**

Let $p_n(z) = Fib_{n+1}(z, -1)$ be a Fibonacci polynomial. These polynomials satisfy a 3-term recurrence with $s(n) = 0$ and $t(n) = 1$. Therefore they are orthogonal with respect to some linear functional $F$.

The corresponding moment sequence is the sequence of aerated Catalan numbers, i.e.

$$F(z^{2n}) = \binom{2n}{n} - \binom{2n}{n-1} = \frac{1}{n+1}\binom{2n}{n} = C_n$$

and $F(z^{2n+1}) = 0$.

Some information about these numbers can be found in the **Online Encyclopedia of Integer Sequences (OEIS) A000108**.

Therefore the Hankel determinants of order 0 and 1 of the aerated Catalan numbers are $d_0(n) = 1$ and $d_1(n) = (-1)^n p_n(0) = (-1)^n Fib_{n+1}(0, -1)$, i.e. $d_1(2n) = (-1)^n$ and $d_1(2n+1) = 0$.

The moments can be computed from the well-known identity (cf. e.g. [6], (1.7))

$$z^n = \sum_{k=0}^{\lfloor \frac{n+1}{2} \rfloor} \left( \binom{n}{k} - \binom{n}{k-1} \right) Fib_{n+1-2k}(z,-1). \quad (1.4)$$

**Example 2**

The normalized Lucas polynomials $p_n(z) = Luc_n(z,-1)$ satisfy a 3-term recurrence with $s(n) = 0$, $t(0) = 2$ and $t(n) = 1$ for $n > 0$. Let $F$ denote the linear functional defined by $F(p_n) = [n = 0]$.

The moment sequence is the sequence of aerated central binomial coefficients (OEIS A126869).

This can be deduced from the identity (cf. e.g. [6], (1.6))

$$z^n = \sum_{k=0}^{\lfloor \frac{n}{2} \rfloor} \binom{n}{k} Luc_{n-2k}(z,-1) \quad (1.5)$$

which gives $F(z^{2n}) = \binom{2n}{n}$ and $F(z^{2n+1}) = 0.$

The corresponding Hankel determinants are $d_0(n) = 2^{n-1}$ and

$$d_1(n) = (-1)^n p_n(0) d_0(n) = (-1)^n Luc_n(0,-1) 2^{n-1},$$

i.e. $d_1(2n) = (-1)^n 4^n$ and $d_1(2n+1) = 0$. For $n = 0$ we set $d_0(0) = 1$.

The generating function of the aerated Catalan numbers is $f(z) = \sum_{n\geq 0} C_n z^{2n} = \dfrac{1-\sqrt{1-4z^2}}{2z^2}$ and satisfies the quadratic equation $f(z) = 1 + z^2 f(z)^2$.

The generating function of the aerated central binomial coefficients is $g(z) = \sum_{n\geq 0} \binom{2n}{n} z^{2n} = \dfrac{1}{\sqrt{1-4z^2}}$

which satisfies $g(z) = 1 + 2z^2 f(z) g(z).$

In this paper we want to study sequences whose generating function satisfies similar equations.

We start with Theorem A which has appeared in the literature in many different disguises (cf. e.g. [1], [8], [10], [14] ). For the sake of completeness we will sketch a proof in Section 2.

**Theorem A**

Let $h(z,x,y) = \sum_{n\geq 0} a(n,x,y)z^n$ satisfy the quadratic equation

$$h(z,x,y) = 1 + (x+y)zh(z,x,y) + xyz^2 h(z,x,y)^2 \tag{1.6}$$

and let more generally

$$H(z,x,y,t) = \sum_{n\geq 0} A(n,x,y,t)z^n = \frac{h(z,x,y)}{1-tzh(z,x,y)}, \tag{1.7}$$

which satisfies

$$H(z,x,y,t) = 1 + (x+y+t)zH(z,x,y,t) + xyz^2 H(z,x,y,t)h(z,x,y) \tag{1.8}$$

and has the explicit expression

$$\begin{aligned}H(z,x,y,t) &= \frac{1-(x+y+2t)z - \sqrt{1-2(x+y)z+(x-y)^2 z^2}}{2z(-t+(t+x)(t+y)z)} \\ &= \frac{2}{1-(x+y+2t)z + \sqrt{1-2(x+y)z+(x-y)^2 z^2}}.\end{aligned} \tag{1.9}$$

Then the Hankel determinants up to order $2$ of the sequence $\bigl(A(n,x,y,t)\bigr)_{n\geq 0}$ are

$$\det\bigl(A(i+j,x,y,t)\bigr)_{i,j=0}^{n-1} = (xy)^{\binom{n}{2}}, \tag{1.10}$$

$$\det\bigl(A(i+j+1,x,y,t)\bigr)_{i,j=0}^{n-1} = (xy)^{\binom{n}{2}}\left(\frac{x^{n+1}-y^{n+1}}{x-y} + t\frac{x^n-y^n}{x-y}\right), \tag{1.11}$$

and

$$\det\bigl(A(i+j+2,x,y,t)\bigr)_{i,j=0}^{n-1} = (xy)^{\binom{n}{2}}\sum_{k=0}^{n}(xy)^{n-k}\left(\frac{x^{k+1}-y^{k+1}}{x-y} + t\frac{x^k-y^k}{x-y}\right)^2. \tag{1.12}$$

**Remark 1.1**

The main example occurs for $(x, y, t) = (1, 1, -1)$. Here $A(n, 1, 1, -1) = C_n$ are the **Catalan numbers**
$C_n = \frac{1}{n+1}\binom{2n}{n}$ (cf. OEIS A000108). Their generating function is

$$H(z,1,1,-1) = C(z) := \sum_{n \geq 0} C_n z^n = \frac{1 - \sqrt{1-4z}}{2z}$$ and satisfies the quadratic equation $C(z) = 1 + zC(z)^2$.

Theorem A generalizes the famous Hankel determinants

$$\det\left(C_{i+j}\right)_{i,j=0}^{n-1} = 1,$$
$$\det\left(C_{i+j+1}\right)_{i,j=0}^{n-1} = 1, \tag{1.13}$$
$$\det\left(C_{i+j+2}\right)_{i,j=0}^{n-1} = n+1.$$

Note that for $x = y$ the expression $\frac{x^n - y^n}{x - y}$ has to be interpreted as $nx^{n-1}$.

It should be noted that in this special case all Hankel determinants are explicitly known (cf. e.g. [10], Theorem 33, [13], [5] or [7]):

$$\det\left(C_{i+j+k}\right)_{i,j=0}^{n-1} = \prod_{j=1}^{k-1}\prod_{i=1}^{j} \frac{2n+j+i}{j+i}. \tag{1.14}$$

But in the general case it is unlikely that such explicit formulae exist. Therefore I restrict myself in this paper to the computation of the Hankel determinants up to order $2$.

An equivalent version of Theorem A is

**Theorem A\***

Let $f(z) = f(z, a, b) = \sum_{n \geq 0} c(n, a, b) z^n$ satisfy the quadratic equation

$$f(z) = 1 + azf(z) + bz^2 f(z)^2, \tag{1.15}$$

*i.e.*

$$f(z) = \frac{1 - az - \sqrt{(1-az)^2 - 4bz^2}}{2bz^2} = \frac{2}{1 - az + \sqrt{(1-az)^2 - 4bz^2}} \tag{1.16}$$

and let

$$\frac{f(z)}{1-tzf(z)} = \frac{2}{1-(a+2t)z+\sqrt{(1-az)^2-4bz^2}} = \sum_{n\geq 0} C(n,a,b,t)z^n. \quad (1.17)$$

*Then the Hankel determinants have the following explicit expression:*

$$\det\left(C(i+j,a,b,t)\right)_{i,j=0}^{n-1} = b^{\binom{n}{2}},$$

$$\det\left(C(i+j+1,a,b,t)\right)_{i,j=0}^{n-1} = b^{\binom{n}{2}}\left(Fib_{n+1}(a,-b)+tFib_n(a,-b)\right), \quad (1.18)$$

$$\det\left(C(i+j+2,a,b,t)\right)_{i,j=0}^{n-1} = b^{\binom{n}{2}}\sum_{j=0}^{n} b^{n-j}\left(Fib_{j+1}(a,-b)+tFib_j(a,-b)\right)^2.$$

*Here $Fib_n(x,s)$ is a **Fibonacci polynomial** defined by the recurrence*
*$Fib_n(x,s) = xFib_{n-1}(x,s) + sFib_{n-2}(x,s)$ with initial values $Fib_0(x,s) = 0$ and $Fib_1(x,s) = 1$.*

To derive this theorem from Theorem A we have only to verify the well-known fact that

$$Fib_n(x+y,-xy) = \frac{x^n - y^n}{x-y}. \quad (1.19)$$

**Remark 1.2**

a) For $(a,b) = (1,0)$ we get as in Example 1 the sequence $(c(n,0,1))_{n\geq 0} = (1,0,1,0,2,0,5,0,14,\cdots)$ of aerated Catalan numbers. The Hankel determinants of order $0$ to $2$ are

$$\det\left(c(i+j,0,1)\right)_{i,j=0}^{n-1} = 1,$$

$$\left(\det\left(c(i+j+1,0,1)\right)_{i,j=0}^{n-1}\right)_{n\geq 0} = (0,1,0,-1,0,1,0,-1,\cdots)$$

$$\left(\det\left(c(i+j+2,0,1)\right)_{i,j=0}^{n-1}\right)_{n\geq 0} = (1,1,2,2,3,3,\cdots).$$

b) For $(x,y) = \left(\frac{1+\sqrt{-3}}{2}, \frac{1-\sqrt{-3}}{2}\right)$ or equivalently $(a,b) = (1,1)$ we get the **Motzkin numbers**
$(M_n) = (1,1,2,4,9,21,51,\cdots)$ (OEIS A001006) with generating function

$$M(z) = \frac{1-z-\sqrt{1-2z-3z^2}}{2z^2} = \frac{2}{1-z+\sqrt{1-2z-3z^2}}. \quad (1.20)$$

Since $\left(Fib_n(1,-1)\right)_{n\geq 0} = (0,1,1,0,-1,-1,0,1,1,0,-1,-1,\cdots)$ is periodic with period 6 the corresponding Hankel determinants are

$$\det\left(M_{i+j}\right)_{i,j=0}^{n-1} = 1,$$
$$\left(\det\left(M_{i+j+1}\right)_{i,j=0}^{n-1}\right)_{n\geq 0} = (1,1,0,-1,-1,0,1,1,0,-1,-1,0,\cdots), \tag{1.21}$$
$$\left(\det\left(M_{i+j+2}\right)_{i,j=0}^{n-1}\right)_{n\geq 0} = (1,2,2,3,4,4,5,6,6,7,8,8,\cdots).$$

An extension of Example 2 is

**Theorem B**

*Let*

$$K(z,x,y) = \sum_{n\geq 0} k(n,x,y)z^n = \frac{1}{\sqrt{(1-(x+y)z)^2 - 4xyz^2}} = \frac{1}{\sqrt{1-2(x+y)z+(x-y)^2 z^2}}.$$

*Then for $n > 0$ the first Hankel determinants are*

$$\det\left(k(i+j,x,y)\right)_{i,j=0}^{n-1} = 2^{n-1}(xy)^{\binom{n}{2}}, \tag{1.22}$$

$$\det\left(k(i+j+1,x,y)\right)_{i,j=0}^{n-1} = \det\left(k(i+j,x,y)\right)_{i,j=0}^{n-1}\left(x^n + y^n\right), \tag{1.23}$$

$$\det\left(k(i+j+2,x,y)\right)_{i,j=0}^{n-1} = \det\left(k(i+j,x,y)\right)_{i,j=0}^{n-1}(xy)^n\left(2 + \sum_{j=1}^{n}\frac{(x^j+y^j)^2}{(xy)^j}\right). \tag{1.24}$$

**Remark 1.3**

The main example is the sequence of **central binomial coefficients** $k(n,1,1) = \binom{2n}{n}$ (OEIS A0001084)

with generating function $K(z,1,1) = \frac{1}{\sqrt{1-4z}} = \sum_{n\geq 0}\binom{2n}{n}z^n$ which we already encountered in Example 2.

The Hankel determinants of orders 0 to 2 are

$$\det\left(\binom{2i+2j}{i+j}\right)_{i,j=0}^{n-1} = 2^{n-1}, \tag{1.25}$$

$$\det\left(\binom{2i+2j+2}{i+j+1}\right)_{i,j=0}^{n-1} = 2^n, \tag{1.26}$$

$$\det\left(\binom{2i+2j+4}{i+j+2}\right)_{i,j=0}^{n-1} = 2^n(2n+1). \tag{1.27}$$

Here too all Hankel determinants are explicitly known (cf. e.g. [5],[7],[10],[13]) :

$$\det\left(\binom{2i+2j+2k}{i+j+k}\right)_{i,j=0}^{n-1} = 2^{n-1+k}\prod_{j=0}^{k-1}\prod_{i=1}^{j}\frac{2n+j+i-1}{i+j}. \tag{1.28}$$

Recall that the **Lucas polynomials** $L_n(x,s)$ are defined by the recursion $L_n(x,s) = xL_{n-1}(x,s) + sL_{n-2}(x,s)$ with initial values $L_0(x,s) = 2$ and $L_1(x,s) = x$ and satisfy $L_n(x+y,-xy) = x^n + y^n$.

In the following we need the **normalized Lucas polynomials** $Luc_n(x,s)$ which coincide with $L_n(x,s)$ for $n > 0$ with initial value $Luc_0(x,s) = 1$. They satisfy $Luc_n(x,s) = xLuc_{n-1}(x,s) + sLuc_{n-2}(x,s)$ for $n > 2$ and $Luc_2(x,s) = xLuc_1(x,s) + 2sLuc_0(x,s)$ with initial values $Luc_0(x,s) = 1$ and $Luc_1(x,s) = x$.

**Theorem B\***

*Let*

$$G(z,a,b) = \sum_{n\geq 0} g(n,a,b)z^n = \frac{1}{\sqrt{(1-az)^2 - 4bz^2}}. \tag{1.29}$$

*Then for $n > 0$ the first Hankel determinants are*

$$\det\left(g(i+j,a,b)\right)_{i,j=0}^{n-1} = 2^{n-1} b^{\binom{n}{2}}, \tag{1.30}$$

$$\det\left(g(i+j+1,a,b)\right)_{i,j=0}^{n-1} = 2^{n-1} b^{\binom{n}{2}} Luc_n(a,-b), \tag{1.31}$$

$$\det\left(g(i+j+2,a,b)\right)_{i,j=0}^{n-1} = 2^{n-1} b^{\binom{n+1}{2}}\left(2 + \sum_{j=1}^{n}\frac{Luc_j(a,-b)^2}{b^j}\right). \tag{1.32}$$

A proof will be given in Section 2.

**Remark 1.4**

As a further example consider the sequence of **central trinomial coefficients** $(g(n,1,1))_{n\geq 0} = (1,1,3,7,19,51,141,\cdots)$ (OEIS A002426). Here we get

$$\det\left(g(i+j,1,1)\right)_{i,j=0}^{n-1} = 2^{n-1}, \tag{1.33}$$

$$\det\left(g(i+j+1,1,1)\right)_{i,j=0}^{n-1} = 2^{n-1} L_n(1,-1), \tag{1.34}$$

where $(L_n(1,-1))_{n\geq 0}$ is the periodic sequence $(2,1,-1,-2,-1,1,\cdots)$ with period $6$, and

$$\det\left(g(i+j+2,1,1)\right)_{i,j=0}^{n-1} = 2^{n-1} r(n) \tag{1.35}$$

with $r(3n) = 6n+2$, $r(3n+1) = 6n+3$, $r(3n+2) = 6n+4$.

It should be noted that we could also use Chebyshev polynomials instead of Fibonacci and Lucas polynomials. The Chebyshev polynomials $T_n(x)$ of the first kind are defined by $T_n(\cos\vartheta) = \cos n\vartheta$ and satisfy $2T_n(x) = L_n(2x,-1)$ and the Chebyshev polynomials of the second kind $U_n(x)$ are defined by $U_n(\cos\vartheta) = \dfrac{\sin(n-1)\vartheta}{\sin\vartheta}$ and satisfy $U_n(x) = Fib_{n+1}(2x,-1)$.

## 2. Hankel determinants and orthogonal polynomials

As already mentioned there is a close connection between Hankel determinants and orthogonal polynomials. This can be sketched in the following way (cf. e.g. [10], [5] or [14]).

Let $(a(n))$ be a sequence of real numbers with $a(0) = 1$.

If all Hankel determinants are unequal to $0$,

$$d_0(n) := \det(a(i+j))_{i,j=0}^{n-1} \neq 0, \tag{2.1}$$

then the monic polynomials

$$p_n(z) := \frac{1}{d_0(n)} \det \begin{pmatrix} a(0) & a(1) & \cdots & a(n-1) & 1 \\ a(1) & a(2) & \cdots & a(n) & z \\ a(2) & a(3) & \cdots & a(n+1) & z^2 \\ \vdots & & & & \vdots \\ a(n) & a(n+1) & \cdots & a(2n-1) & z^n \end{pmatrix} \tag{2.2}$$

are orthogonal with respect to the linear functional $F$ defined by $F(z^n) = a(n)$. Note that we define $d_0(0) = 1$.

By Favard's theorem there exist numbers $s(n), t(n)$ such that

$$p_n(z) = (z - s(n-1))p_{n-1}(z) - t(n-2)p_{n-2}(z). \tag{2.3}$$

The sequence $(a(n))$ is uniquely determined by the sequences $s(n)$ and $t(n)$. For define $v(n, j)$ as the uniquely determined coefficients in

$$z^n = \sum_{j=0}^{n} v(n, j) p(j, z). \tag{2.4}$$

Comparing coefficients in $z^n = zz^{n-1} = \sum_{j=0}^{n} v(n-1, j)zp(j, z)$ this is equivalent with

$$\begin{aligned} v(0, j) &= [j = 0] \\ v(n, 0) &= s(0)v(n-1, 0) + t(0)v(n-1, 1) \\ v(n, j) &= v(n-1, j-1) + s(j)v(n-1, j) + t(j)v(n-1, j+1). \end{aligned} \tag{2.5}$$

The numbers $v(n, j)$ can be interpreted as weights of some lattice paths: Consider non-negative lattice paths from $(0,0)$ to $(n, j)$ with horizontal steps $H = (1,0)$, up-steps $U = (1,1)$ and down-steps $D = (1,-1)$. Define the weight $w_0$ of a path as the product of all steps of the path, where $w_0(U) = 1$, $w_0(H) = s(k)$ if $H$ lies on height $k$ and $w_0(D) = t(k)$ if the endpoint of $D$ in on height $k$. The weight of a set of paths is defined as the sum of their weights. Then $v(n, j)$ is the weight of the set of all lattice paths from $(0,0)$ to $(n, j)$ and $v(n, 0) = a(n)$.

The Hankel determinant $d_0(n) = \det\left(v(i+j, 0)\right)_{i,j=0}^{n-1}$ is given by

$$d_0(n) = \det\left(a(i+j)\right)_{i,j=0}^{n-1} = \prod_{i=1}^{n-1} \prod_{j=0}^{i-1} t(j). \tag{2.6}$$

To prove this observe that each non-negative lattice path from $(0,0)$ to $(m+n, 0)$ must contain a unique point $(m, i)$ with $0 \leq i \leq m$. The weight of all paths from $(0,0)$ to $(m, i)$ is $v(m, i)$ and the weight of all paths from $(m, i)$ to $(m+n, 0)$ is the same as the weight $v_1(n, i)$ of all paths from $(0,0)$ to $(n, i)$ with the weights $w_1(D) = 1$, $w_1(U) = t(k)$ if the initial point of $U$ is on height $k$ and $w_1(H) = s(k)$ if $H$ is on height $k$. Then $v_1(n, i) = v(n, i) \prod_{j=0}^{i-1} t(j)$ due to $i$ upsteps.

This implies $v(m+n,0) = \sum_{i=0}^{m} v(m,i)v(n,i)\prod_{j=0}^{i-1} t(j)$.

If we denote by $A_n$ the lower triangular matrix $A_n = (v(i,j))_{i,j=0}^{n-1}$ with diagonal $v(i,i) = 1$ and by $B_n$ the diagonal matrix with entries $b(i,i) = \prod_{j=0}^{i-1} t(j)$ then $(v(i+j,0))_{i,j=0}^{n-1} = A_n B_n A_n^t$. Taking determinants we get (2.6).

From (2.2) it is clear that

$$d_1(n) = \det(a(i+j+1))_{i,j=0}^{n-1} = (-1)^n p_n(0) d_0(n). \tag{2.7}$$

From the **condensation formula for determinants** (cf. [10], p.17) we get for the Hankel determinants $d_2(n) = \det(a(i+j+2))_{i,j=0}^{n-1}$

$$d_2(n) d_0(n) = d_0(n+1) d_2(n-1) + d_1(n)^2. \tag{2.8}$$

**Proof of Theorem A**

If we set

$$\sum_{n \geq 0} v(n,k) z^n := z^k H(z,x,y,t) h(z,x,y)^k \tag{2.9}$$

then (2.5) is satisfied with

$$s(0) = x + y + t, s(n) = x + y, t(n) = xy. \tag{2.10}$$

This is seen by comparing coefficients in (1.8) and in (1.6) after multiplication by $z^k H(z,x,y,t) h(z,x,y)^{k-1}$.

Thus we get (1.10).

Moreover we see that

$$p_n(z) = Fib_{n+1}(z-(x+y),-xy) - t Fib_n(z-(x+y),-xy). \tag{2.11}$$

Therefore $r(n) = (-1)^n p_n(0)$ satisfies the recurrence $r(n) = (x+y)r(n-1) - (xy)r(n-2)$ with initial values $r(0) = 1 = \frac{x-y}{x-y} + t\frac{x^0-y^0}{x-y}$ and $r(1) = x+y+t = \frac{x^2-y^2}{x-y} + t\frac{x-y}{x-y}$ and since

$\frac{x^{n+1}-y^{n+1}}{x-y} = (x+y)\frac{x^n-y^n}{x-y} - xy\frac{x^{n-1}-y^{n-1}}{x-y}$ we see that $r(n) = \frac{x^{n+1}-y^{n+1}}{x-y} + t\frac{x^n-y^n}{x-y}$.

This implies (1.11). To show (1.12) we use the condensation formula and obtain

$$d_2^{(2)}(n, x+y, xy) = (xy)^n d_2^{(2)}(n-1, x+y, xy) + (xy)^{\binom{n}{2}} d_1^{(2)}(n, x+y, xy)^2.$$

With induction we get $d_2^{(2)}(n, x+y, xy) = (xy)^{\binom{n}{2}} \sum_{k=0}^{n} (xy)^{n-k} \left( \frac{x^{k+1} - y^{k+1}}{x-y} + t\frac{x^k - y^k}{x-y} \right)^2.$

**Remark 2.1**

For $x = i, y = -i, t = 0$ we get the sequence of aerated Catalan numbers $(1,0,1,0,2,0,5,0,\cdots)$. In this case we get (cf. e.g. [5],[6],[7] or OEIS A053121)

$$v(2n, 2j) = \binom{2n}{n-j} - \binom{2n}{n-j-1}, v(2n+1, 2j+1) = \binom{2n+1}{n-j} - \binom{2n+1}{n-j-1} \text{ and } v(n,j) = 0 \text{ if}$$
$(n-j) \equiv 1 \mod 2.$

The corresponding orthogonal polynomials are $p_n(z) = Fib_{n+1}(z,-1) = U_n\left(\frac{z}{2}\right).$

A comparison with the orthogonality relation for the Chebyshev polynomials of the second kind

$$\frac{2}{\pi} \int_{-1}^{1} U_n(x) U_m(x) \sqrt{1-x^2} \, dx = [n=m] \text{ gives}$$

$$\frac{1}{2\pi} \int_{-2}^{2} F_{n+1}(x,-1) \sqrt{4-x^2} \, dx = [n=0]. \quad (2.12)$$

Thus we get a representation of the linear functional $F$ as an integral

$$F(p) = \frac{1}{2\pi} \int_{-2}^{2} p(x) \sqrt{4-x^2} \, dx \quad (2.13)$$

for each polynomial $p(x)$.

**Remark 2.2**

For $a(n) = c(n,a,b)$ the orthogonal polynomials corresponding to Theorem A* satisfy

$p_n(z) = (z-a) p_{n-1}(z) - b p_{n-2}(z)$ with initial values $p_{-1}(z) = 0$ and $p_0(z) = 1$, i.e.
$p_n(z) = Fib_{n+1}(z-a, -b).$

Thus

$$\frac{1}{b^{\binom{n}{2}}}\det\begin{pmatrix} c(0,a,b) & c(1,a,b) & \cdots & c(n-1,a,b) & 1 \\ c(1,a,b) & c(2,a,b) & \cdots & c(n,a,b) & z \\ c(2,a,b) & c(3,a,b) & \cdots & c(n+1,a,b) & z^2 \\ \vdots & & & & \vdots \\ c(n,a,b) & c(n+1,a,b) & \cdots & c(2n-1,a,b) & z^n \end{pmatrix} = Fib_{n+1}(z-a,-b). \quad (2.14)$$

By (2.12) we get

$$c(n,a,b) = \frac{1}{2\pi b}\int_{a-2\sqrt{b}}^{a+2\sqrt{b}} x^n \sqrt{4b-(x-a)^2}\,dx = \frac{1}{2\pi}\int_{-2}^{2}(a+z\sqrt{b})^n\sqrt{4-z^2}\,dz. \quad (2.15)$$

For $a(n) = C(n,a,b,t)$ we get in the same way $p_n(z) = Fib_{n+1}(z-a,-b) - tFib_n(z-a,-b)$.

In the special case $(a,b,t) = (2,1,-1)$ we get $C(n,2,1,-1) = C_n$ and

$$p_n(z) = Fib_{n+1}(z-2,-1) + Fib_n(z-2,-1) = Fib_{2n+1}(\sqrt{z},-1).$$

**Proof of Theorem B\***

To prove Theorem B* we define $v(n,j)$ by the generating functions

$$V_j(z) = \sum_{n\geq 0} v(n,j)z^n = z^j G(z,a,b) f(z,a,b)^j. \quad (2.16)$$

Then

$$\begin{aligned} V_0(z) &= 1 + azV_0(z) + 2bzV_1(z), \\ V_j(z) &= z\left(V_{j-1}(z) + aV_j(z) + bV_{j+1}(z)\right). \end{aligned} \quad (2.17)$$

The first line means that $G(z,a,b) = 1 + azG(z,a,b) + 2bz^2 G(z,a,b)f(z,a,b)$

or equivalently $1 - az - 2bz^2 \dfrac{1-az-\sqrt{(1-az)^2-4bz^2}}{2bz^2} = \sqrt{(1-az)^2-4bz^2}$,

which is obvious. The second line is simply the quadratic equation satisfied by $f(z,a,b)$.

Therefore (2.5) holds with

$$s(n) = a, t(0) = 2b \quad (2.18)$$

and

$$t(n) = b \quad (2.19)$$

for $n > 0$.

This means that the corresponding normalized orthogonal polynomials are

$$p_n(z) = Luc_n(z-a,-b). \tag{2.20}$$

By (2.6) we get (1.30).

From (2.7) we deduce (1.31), because $r(n) = (-1)^n p(n,0)$ satisfies $r(n) = ar(n-1) - br(n-2)$ with initial values $r(1) = a$ and $r(2) = a^2 - 2b$. Since $L_n(a,-b)$ satisfies the same recurrence and coincides with $r(n)$ for $n=1$ and $n=2$ we see that $r(n) = Luc_n(a,-b)$.

Identity (1.32) is a simple consequence of the condensation formula.

**Remark 2.3**

For $(a,b) = (2,1)$ it is easily verified (cf. e.g. [5],[7]) that (2.5) holds with $v(n,j) = \binom{2n}{n-j}$. In this case we have $f(z,2,1) = c(z) = \dfrac{1-2z-\sqrt{1-4z}}{2z^2} = \dfrac{C(z)-1}{z} = \sum_{n \geq 0} C_{n+1} z^n$.

Comparing with (2.16) we get the well-known formula

$$\sum_{n \geq 0} \binom{2n}{n-j} z^n = \frac{z^j c(z)^j}{\sqrt{1-4z}} = \frac{(C(z)-1)^j}{\sqrt{1-4z}}. \tag{2.21}$$

**Remark 2.4**

The orthogonal polynomials belonging to Theorem B* are $p_n(z) = Luc_n(z-a,-b)$. Thus for $n > 0$

$$\frac{1}{2^{n-1} b^{\binom{n}{2}}} \det \begin{pmatrix} g(0,a,b) & g(1,a,b) & \cdots & g(n-1,a,b) & 1 \\ g(1,a,b) & g(2,a,b) & \cdots & g(n,a,b) & z \\ g(2,a,b) & g(3,a,b) & \cdots & g(n+1,a,b) & z^2 \\ \vdots & & & & \vdots \\ g(n,a,b) & g(n+1,a,b) & \cdots & g(2n-1,a,b) & z^n \end{pmatrix} = Luc_n(z-a,-b). \tag{2.22}$$

For $(a,b) = (2,1)$ we get $g(n,2,1) = \binom{2n}{n}$ and $p_n(z) = Luc_n(z-2,-1) = Luc_{2n}(\sqrt{z},-1)$.

For $(a,b) = (0,1)$ we get Example 2: $g(2n,0,1) = \binom{2n}{n}$ and $g(2n+1,0,1) = 0$ and $p_n(z) = Luc_n(z,-1)$.

Now observe that the Chebyshev polynomials $T_n(\cos\vartheta) = \cos n\vartheta$ satisfy $2T_n(x) = L_n(2x,-1)$ and

$$\frac{1}{\pi}\int_{-1}^{1}\frac{T_n(x)}{\sqrt{1-x^2}}dx = [n=0].$$

This implies that

$$F(Luc_n(z,-1)) = \frac{1}{\pi}\int_{-2}^{2}\frac{Luc_n(x,-1)}{\sqrt{4-x^2}}dx = [n=0] \tag{2.23}$$

and therefore

$$F(z^{2n}) = \frac{1}{\pi}\int_{-2}^{2}\frac{x^{2n}}{\sqrt{4-x^2}}dx = \binom{2n}{n} \text{ and } F(z^{2n+1}) = \frac{1}{\pi}\int_{-2}^{2}\frac{x^{2n+1}}{\sqrt{4-x^2}}dx = 0.$$

More generally we have

$$g(n,a,b) = \frac{1}{\pi}\int_{-2}^{2}\frac{(a+z\sqrt{b})^n}{\sqrt{4-z^2}}dz. \tag{2.24}$$

### 3. Some general observations

Let $m$ be a positive integer. Our main purpose is to study Hankel determinants of sequences $(c(n,m,a,b))$ and $(C(n,m,a,b,t))$ respectively whose generating functions satisfy an equation of the form

$$f_m(z,a,b) = \sum_{n\geq}c(n,m,a,b)z^n = 1 + azf_m(z,a,b) + bz^m f_m(z,a,b)^2 \tag{3.1}$$

and

$$F_m(z,a,b,t) = \sum_{n\geq 0}C(n,m,a,b,t)z^n = \frac{f_m(z,a,b)}{1-tzf_m(z,a,b)}. \tag{3.2}$$

The coefficients $C(n,m,a,b,t)$ can be regarded as generating functions of weighted non-negative lattice paths:

In the sequel for given $m$ we always consider lattice paths from $(0,0)$ to $(n,0)$ with horizontal steps $H = (1,0)$, up-steps $U = (1,1)$ and down-steps $D = (m-1,-1)$ of width $m-1$.

Define the **weight** $w_0$ of a path as the product of all steps of the path, where $w_0(U) = 1$, $w_0(H) = s(k)$ if $H$ lies on height $k$ and $w_0(D) = t(k)$ if the endpoint of $D$ in on height $k$.

The weight of a set of paths is defined as the sum of their weights.

**Lemma 3.1**

Let $f(z)$ be the generating function of all non-negative paths which start and end on the $x-axis$ and let $g(z)$ be the generating function of all paths which start at height $1$ and end on height $1$ and never reach height $0$. Then

$$f(z) = 1 + s(0)zf(z) + t(0)z^m f(z)g(z). \qquad (3.3)$$

**Proof**

This is clear since the generating function of all non-trivial paths which start with a horizontal step is $s(0)zf(z)$ and each other path has a unique decomposition of the following form: an up-step, a maximal path which never falls under height 1 and ends on height 1, a down-step to height 0 and a further non-negative path which ends on height 0.

**Lemma 3.2**

Let $s(n) = a$ and $t(n) = b$ and let $c(n,m,a,b)$ be the weight of all non-negative paths from $(0,0)$ to $(n,0)$ and let $f_m(z,a,b) = \sum_{n \geq 0} c(n,m,a,b)z^n$ be the generating function of the weights.

Then for all positive integers $m$

$$f_m(z,a,b) = 1 + azf_m(z,a,b) + bz^m f_m(z,a,b)^2. \qquad (3.4)$$

The weights $c(n,m,a,b)$ are given by

$$c(n,m,a,b) = \sum_{k=0}^{\frac{n}{m}} \binom{2k}{k} \frac{1}{k+1} \binom{n+(2-m)k}{2k} a^{n-mk} b^k. \qquad (3.5)$$

**Proof**

There are $\binom{2k}{k}\frac{1}{k+1}$ different possible sets of $k$ up-steps and $k$ downsteps.

For each fixed set of $k$ up-steps and $k$ downsteps there are $\binom{(n-mk)+2k}{2k}$ different positions for the remaining $n-mk$ horizontal steps.

Note that for $m = 2$ (3.5) is the same as (2.15).

Let somewhat more generally

$$F_m(z,a,b,t) = \sum_{n \geq 0} C(n,m,a,b,t)z^n = \frac{f_m(z,a,b)}{1-tzf_m(z,a,b)}. \tag{3.6}$$

A closed formula is

$$F_m(z,a,b,t) = \frac{2}{1-(a+2t)z + \sqrt{(1-az)^2 - 4bz^m}}. \tag{3.7}$$

Then it is easily verified that

$$F_m(z,a,b,t) = 1 + (a+t)zF_m(z,a,b,t) + bz^m F_m(z,a,b,t) f_m(z,a,b). \tag{3.8}$$

We are interested in the Hankel determinants $d_k^{(m)}(n,a,b) = \det\left(c(i+j+k,m,a,b)\right)_{i,j=0}^{n-1}$ and slightly more generally in $D_k^{(m)}(n,a,b,t) = \det\left(C(i+j+k,m,a,b,t)\right)_{i,j=0}^{n-1}$.

**Unrestricted paths**

A related question concerns the case of all paths.

Let $g(n,m,a,b)$ be the weight of all paths from $(0,0)$ to $(n,0)$. If there are $k$ up-steps there must also be $k$ down-steps and $n-km$ horizontal steps. The number of steps is therefore $n+2k-mk$. There are $\frac{(n+2k-mk)!}{k!k!(n-mk)!} = \binom{2k}{k}\binom{n+2k-mk}{2k}$ different possibilities for the distribution of the steps. Thus

$$g(n,m,a,b) = \sum_{k=0}^{\lfloor \frac{n}{m} \rfloor} \binom{2k}{k}\binom{n+2k-mk}{2k} a^{n-mk}b^k. \tag{3.9}$$

**Remark 3.3**

More information about such paths can be found in [2]. For $m=2$ (3.9) is the same as (2.24).

It should be noted that there are also other paths which give the same weights. For example consider all paths from $(0,0)$ to $(n,n)$ with horizontal steps $(m,0)$, vertical steps $(0,m)$ and diagonal steps $(1,1)$. If there are $k$ horizontal steps there must also be $k$ vertical steps and $n-mk$ diagonal steps. Thus we get the same formula.

The generating function can be obtained in the same manner as in [11], Exercise 5.b.

$$G_m(z,a,b) = \sum_{n \geq 0} g(n,m,a,b)z^n = \sum_{n \geq 0} z^n \sum_{k=0}^{\lfloor n/m \rfloor} \binom{2k}{k}\binom{n+2k-mk}{2k} a^{n-mk} b^k$$

$$= \sum_{k \geq 0} \binom{2k}{k} \sum_{n \geq km} \binom{n+2k-mk}{2k} a^{n-mk} b^k z^n = \sum_{k \geq 0} \binom{2k}{k}\left(\frac{b}{a^m}\right)^k \sum_{n \geq km} \binom{n+2k-mk}{2k} a^n z^n$$

$$= \sum_{k \geq 0} \binom{2k}{k}\left(\frac{b}{a^m}\right)^k \sum_{j \geq 0} \binom{2k+j}{2k} a^{j+mk} z^{j+mk} = \sum_{k \geq 0} \binom{2k}{k}\left(\frac{bz^{m-2}}{a^2}\right)^k \sum_{j \geq 0} \binom{2k+j}{2k}(az)^{2k+j}$$

$$= \sum_{k \geq 0} \binom{2k}{k}\left(\frac{bz^{m-2}}{a^2}\right)^k \frac{(az)^{2k}}{(1-az)^{2k+1}} = \frac{1}{1-az} \sum_{k \geq 0} \binom{2k}{k}\left(\frac{bz^m}{(1-az)^2}\right)^k = \frac{1}{\sqrt{(1-az)^2 - 4bz^m}}.$$

Here too we are interested in computing the Hankel determinants

$$dd_k^{(m)}(n,a,b) = \det\left(g(i+j+k,m,a,b)\right)_{i,j=0}^{n-1}. \tag{3.10}$$

We found some facts for $m > 2$ but could not find proofs for our conjectures.

## 4. The special case $m=1$.

The case $m=1$ can be reduced to $m=2$ since by (3.7) it is easily seen that

$$f_1(z,a,b) = F_1(z,a,b,0) = F_2(z, a+2b, b(a+b), -b). \tag{4.1}$$

Therefore we have $s(0) = a+b$, $s(n) = a+2b$ for $n > 0$ and $t(n) = b(a+b)$.

This implies

**Theorem 4.1**
*The Hankel determinants $d_k^{(1)}(n,a,b)$ are*

$$d_0^{(1)}(n,a,b) = b^{\binom{n}{2}}(a+b)^{\binom{n}{2}},$$
$$d_1^{(1)}(n,a,b) = b^{\binom{n}{2}}(a+b)^{\binom{n+1}{2}}, \tag{4.2}$$
$$d_2^{(1)}(n,a,b) = b^{\binom{n}{2}}(a+b)^{\binom{n+1}{2}} \frac{(a+b)^{n+1} - b^{n+1}}{a}.$$

**Remark 4.2**

The main example is the sequence of **large Schröder numbers** $(R_n)_{n \geq 0} = (1, 2, 6, 22, 90, \cdots)$ (OEIS A006318) with $R_n = c(n,1,1,1)$. Their generating function satisfies $F(z) = 1 + zF(z) + zF(z)^2$.

The corresponding Hankel determinants are (cf. e.g.[4] )

$$\det\left(R_{i+j}\right)_{i,j=0}^{n-1} = 2^{\binom{n}{2}}, \tag{4.3}$$

$$\det\left(R_{i+j+1}\right)_{i,j=0}^{n-1} = 2^{\binom{n+1}{2}}, \tag{4.4}$$

$$\det\left(R_{i+j+2}\right)_{i,j=0}^{n-1} = 2^{\binom{n+1}{2}}\left(2^{n+1}-1\right), \tag{4.5}$$

More generally we get

**Corollary 4.3**

$$D_0^{(1)}(n,a,b,t) = b^{\binom{n}{2}}(a+b)^{\binom{n}{2}},$$

$$D_1^{(1)}(n,a,b,t) = b^{\binom{n}{2}}(a+b)^{\binom{n}{2}}\left((a+b)^n + t\frac{(a+b)^n - b^n}{a}\right), \tag{4.6}$$

$$D_2^{(1)}(n,a,b,t) = b^{\binom{n}{2}}(a+b)^{\binom{n}{2}}\sum_{j=0}^{n}(a+b)^{n-j}b^{n-j}\left((a+b)^j + t\frac{(a+b)^j - b^j}{a}\right)^2.$$

**Proof.**

By (1.19) we have $Fib_n(a+2b, -b(a+b)) = \dfrac{(a+b)^n - b^n}{a}$ and therefore
$Fib_{n+1}(a+2b, -b(a+b)) - bFib_n(a+2b, -b(a+b)) = (a+b)^n$.

Further we have

$$\sum_{j=0}^{n}(b(a+b))^{\binom{n+1}{2}-j}\left(Fib_{j+1}(a+2b,-b(a+b)) - bFib_j(a+2b,-b(a+b))\right)^2$$

$$= b^{\binom{n}{2}}(a+b)^{\binom{n+1}{2}}\sum_{j=0}^{n}b^{n-j}(a+b)^j = b^{\binom{n}{2}}(a+b)^{\binom{n+1}{2}}\frac{(a+b)^{n+1} - b^{n+1}}{a}.$$

In the same way we can compute the determinants $dd_k^{(1)}(n,a,b)$ since

$$G_1(z,a,b) = \sum_{n\geq 0} g(n,1,a,b)z^n = \frac{1}{\sqrt{(1-az)^2 - 4bz}}$$

$$= \frac{1}{\sqrt{(1-(a+2b)z)^2 - 4b(a+b)z^2}} = G_2(z, a+2b, b(a+b)).$$

**Theorem 4.4**

$$dd_0^{(1)}(n,a,b) = 2^{n-1} b^{\binom{n}{2}} (a+b)^{\binom{n}{2}}, \tag{4.7}$$

$$dd_1^{(1)}(n,a,b) = 2^{n-1} b^{\binom{n}{2}} (a+b)^{\binom{n}{2}} \left( (a+b)^n + b^n \right), \tag{4.8}$$

$$\begin{aligned} dd_2^{(1)}(n,0,a,b) &= dd_0^{(1)}(n,0,a,b) \left( 2b^n (a+b)^n + \sum_{j=1}^{n} \left( (a+b)^j + b^j \right)^2 b^{n-j} (a+b)^{n-j} \right) \\ &= dd_0^{(1)}(n,0,a,b) \left( \frac{(a+b)^{2n+1} - b^{2n+1}}{a} + (2n+1) b^n (a+b)^n \right) \end{aligned} \tag{4.9}$$

*for* $n > 0$.

Note that $L_n(x+y, -xy) = x^n + y^n$.

**Remark 4.5**

The numbers $\left( g(n,1,1) \right)_{n \geq 0} = (1, 3, 13, 63, 321, \cdots)$ (OEIS A001850) are the **central Delannoy numbers**.

Their Hankel determinants are

$$\begin{aligned} dd_0^{(1)}(n,1,1) &= 2^{\binom{n+1}{2} - 1}, \\ dd_1^{(1)}(n,1,1) &= 2^{\binom{n+1}{2} - 1} \left( (2^n + 1) \right), \\ dd_2^{(1)}(n,1,1) &= 2^{\binom{n+1}{2} - 1} \left( 2^{2n+1} - 1 + (2n+1) 2^n \right). \end{aligned} \tag{4.10}$$

**Remark 4.6.**

For $g(n,1,0,b) = \binom{2n}{n} b^n$ we get the Hankel determinants

$$\begin{aligned} dd_0^{(1)}(n,0,b) &= 2^{n-1} b^{n(n-1)}, \\ dd_1^{(1)}(n,0,b) &= 2^n b^{n^2}, \\ dd_2^{(1)}(n,0,b) &= 2^n (2n+1) b^{n^2+n}. \end{aligned} \tag{4.11}$$

## 5. Some observations for $m > 2$.

Now we want to study the case $m > 2$. Here some of the determinants $d_0^{(m)}(n,a,b)$ and $dd_0^{(m)}(n,a,b)$ vanish.

More precisely we have

**Theorem 5.1**

$$d_0^{(m)}(mn+k,a,b) = 0 \text{ for } 2 \leq k \leq m-1,$$

$$d_1^{(m)}(mn+k,a,b) = 0 \text{ for } 2 \leq k \leq m-2,$$

$$d_2^{(m)}(mn+k,a,b) = 0 \text{ for } 2 \leq k \leq m-3.$$

*The same holds for the determinants $dd_j^{(m)}(mn+k,a,b)$.*

**Proof**

It suffices to show that there are polynomials $h(n,m,j)$ in the variables $a,b$ such that for $1 \leq i \leq mn+m-1$

$$\sum_{j=0}^{mn+1} h(n,m,j) c(nm+i-j,m,a,b) = 0 \tag{5.1}$$

or equivalently that there exist polynomials $P(n,m,z)$ of degree $mn+1$ in $z$ such that for $1 \leq i \leq mn+m-1$

$$[z^{nm+i}] P(n,m,z) f_m(z,a,b) = 0. \tag{5.2}$$

For then we have $d_0^{(m)}(mn+k,a,b) = 0$ for $2 \leq k \leq m-1$.

For $2 \leq k \leq m-2$ this implies also that $d_1^{(m)}(mn+k,a,b) = 0$ and for $2 \leq k \leq m-3$ that $d_2^{(m)}(mn+k,a,b) = 0$.

To show this we prove the following

**Lemma 5.2**

*Let*

$$P(n,m,z) = Fib_{2n+2}(1-az,-bz^m) = \sum_{j=0}^{n} \binom{2n+1-j}{j} (1-az)^{2n+1-2j} (-bz^m)^j. \tag{5.3}$$

*Then $P(n,m,z)$ is a polynomial in $z$ with degree $nm+1$ which satisfies*

$$[z^{nm+k}] P(n,m,z) f_m(z,a,b) = 0 \tag{5.4}$$

*for $1 \leq k \leq mn+m-1$.*

**Proof of Lemma 5.2**

The degree of $P(n,m,z)$ is $\max_{0 \le j \le n}(2n+1-2j+mj) = mn+1$.

**Binet's formula** for the Fibonacci polynomials $Fib_n(x,s) = \dfrac{\left(\dfrac{x+\sqrt{x^2+4s}}{2}\right)^n - \left(\dfrac{x-\sqrt{x^2+4s}}{2}\right)^n}{\sqrt{x^2+4s}}$ gives

$$Fib_{2n+2}(1-az,-bz^m) = \frac{\alpha(m,z)^{2n+2} - \beta(m,z)^{2n+2}}{\alpha(m,z) - \beta(m,z)} \tag{5.5}$$

with

$$\alpha(m,z) = \alpha(m,z,a,b) = \frac{1-az + \sqrt{(1-az)^2 - 4bz^m}}{2} \tag{5.6}$$

and

$$\beta(m,z) = \beta(m,z,a,b) = \frac{1-az - \sqrt{(1-az)^2 - 4bz^m}}{2}. \tag{5.7}$$

Observe that $\alpha(m,z)$ and $f_m(z,a,b) = \dfrac{1}{\alpha(m,z)}$ are power series with constant term 1 and that $\beta(m,z)$ is a power series with smallest power $m$ since $\alpha(m,z)\beta(m,z) = bz^m$.

This implies that

$$P(n,m,z) f_m(z,a,b) = \frac{\alpha(m,z)^{2n+2} - \beta(m,z)^{2n+2}}{\alpha(m,z) - \beta(m,z)} \frac{1}{\alpha(m,z)}$$

$$= \frac{1}{\alpha(m,z) - \beta(m,z)} \left( \alpha(m,z)^{2n+1} - \frac{\beta(m,z)^{2n+2}}{\alpha(m,z)} \right)$$

$$= \frac{\alpha(m,z)^{2n+1} - \beta(m,z)^{2n+1}}{\alpha(m,z) - \beta(m,z)} + \frac{\beta(m,z)^{2n+1}}{\sqrt{(1-az)^2 - 4bz^m}} \left( 1 - \frac{\beta(m,z)}{\alpha(m,z)} \right)$$

$$= Fib_{2n+1}(1-az,-bz^m) + z^{(2n+1)m} r(z)$$

for some formal power series $r(z)$.

Since $\deg Fib_{2n+1}(1-az,-bz^m) = mn$ we see that for $1 \le k \le mn+m-1$

$[z^{nm+k}] P(n,m,z) f_m(z,a,b) = 0.$

**Remark 5.3**

For $F_m(z,a,b,t)$ instead of $P(n,m,z)$ we have to choose
$$P(n,m,z,t) = Fib_{2n+2}(1-az,-bz^m) - tzFib_{2n+1}(1-az,-bz^m).$$

For in this case $F_m(z,a,b,t) = \dfrac{1}{\alpha(m,z)-tz}$ and therefore

$$P(n,m,z,t)f_m(z,a,b) = \frac{(\alpha(m,z)-tz)\alpha(m,z)^{2n+1} - (\beta(m,z)-tz)\beta(m,z)^{2n+1}}{\alpha(m,z)-\beta(m,z)} \frac{1}{\alpha(m,z)-tz}$$
$$= Fib_{2n+1}(1-az,-bz^m) + z^{(2n+1)m}r(z).$$

To prove the assertion for the determinants $dd_j^{(m)}(mn+k,a,b)$ observe that Binet's formula for the Lucas polynomials gives $L_n(x,s) = \left(\dfrac{x+\sqrt{x^2+4s}}{2}\right)^n + \left(\dfrac{x-\sqrt{x^2+4s}}{2}\right)^n$.

We claim that

$$[z^{nm+k}]L_{2n+1}(1-az,-bz^m)\frac{1}{\sqrt{(1-az)^2-4bz^m}} = 0 \text{ for } 1 \le k \le mn+m-1.$$

But this follows immediately from the fact that

$$L_{2n+1}(1-az,-bz^m)\frac{1}{\sqrt{(1-az)^2-4bz^m}} = \frac{\alpha(m,z)^{2n+1}+\beta(m,z)^{2n+1}}{\alpha(m,z)-\beta(m,z)}$$
$$= \frac{\alpha(m,z)^{2n+1}-\beta(m,z)^{2n+1}}{\alpha(m,z)-\beta(m,z)} + 2\frac{\beta(m,z)^{2n+1}}{\sqrt{(1-az)^2-4bz^m}} = Fib_{2n+1}(1-az,-bz^m) + z^{(2n+1)m}r(z)$$

for some power series $r(z)$.

Since some of the determinants $d_0^{(m)}(n,a,b)$ vanish the method of orthogonal polynomials is not applicable.

In some cases elementary matrix manipulations will do the task. Consider e.g. the determinants
$$d_0^{(m)}(n,0,b) = \det\bigl(c(i+j,m,0,b)\bigr)_{i,j=0}^{n-1}.$$

Here we have $f_m(z,0,b) = \sum\limits_{n\ge 0} c(n,m,0,b)z^n = 1 + bz^m f_m(z,0,b)^2$ and therefore

$c(nm,m,0,b) = C_n b^n$ and $c(n,m,0,b) = 0$ if $n \not\equiv 0 \bmod m$.

**Theorem 5.4**

*Let $m \geq 1$. Then*

$$d_0^{(m)}(mn,0,b) = (-1)^{\binom{m-1}{2}n} b^{n(mn-1)}$$

$$d_0^{(m)}(mn+1,0,b) = (-1)^{\binom{m-1}{2}n} b^{n(mn+1)}$$

(5.8)

*and $d_0^{(m)}(n,0,b) = 0$ for all other $n$.*

**Proof**

Because of Theorem 5.1 we need only show (5.8).

For $m = 1$ we already know that $d_0^{(1)}(n,0,b) = b^{n^2-n}$ and that $d_1^{(1)}(n,0,b) = b^{n^2}$.

Consider first the matrix $\left(c(i+j,m,0,b)\right)_{i,j=0}^{mn-1}$. We write its rows $r_i$ in the order
$r_0, r_m, \cdots, r_{(n-1)m}, r_1, r_{m+1}, \cdots, r_{(n-1)m+1}, \cdots, r_{m-1}, r_{2m-1}, \cdots, r_{(n-1)m+m-1}$. Then we change the colums $s_j$ of the new matrix in the order $s_0, s_m, \cdots, s_{(n-1)m}, s_{m-1}, s_{2m-1}, \cdots, s_{(n-1)m+m-1}, \cdots, s_1, s_{m+1}, \cdots, s_{(n-1)m+1}$.

The new matrix is of the form

$$\begin{pmatrix} A_0 & & & \\ & A_1 & & \\ & & \ddots & \\ & & & A_{m-1} \end{pmatrix},$$

where $A_0 = \left(C_{i+j} b^{i+j}\right)_{i,j=0}^{n-1}$ and the other $A_\ell$ are equal to $A_\ell = \left(C_{i+j+1} b^{i+j+1}\right)_{i,j=0}^{n-1}$.

Thus $d_0^{(m)}(mn,0,b) = (-1)^{n\binom{m-1}{2}} \prod_{j=0}^{m-1} \det A_j = (-1)^{n\binom{m-1}{2}} b^{n(n-1)} \left(b^{n^2}\right)^{m-1} = (-1)^{n\binom{m-1}{2}} b^{n(nm-1)}.$

For the matrix $\left(c(i+j,m,0,b)\right)_{i,j=0}^{mn}$ we have an analogous decomposition. The only difference is that $A_0 = \left(C_{i+j} b^{i+j}\right)_{i,j=0}^{n}$. This implies that

$$d_0^{(m)}(mn+1,0,b) = (-1)^{n\binom{m-1}{2}} \prod_{j=0}^{m-1} \det A_j = (-1)^{n\binom{m-1}{2}} b^{n(n+1)} \left(b^{n^2}\right)^{m-1} = (-1)^{n\binom{m-1}{2}} b^{n(nm+1)}.$$

In the same way using (4.11) we get

**Theorem 5.5**

*Let $m \geq 1$. Then*

$$dd_0^{(m)}(mn,0,b) = (-1)^{\binom{m-1}{2}n} 2^{mn-1} b^{n(mn-1)}$$

$$dd_0^{(m)}(mn+1,0,b) = (-1)^{\binom{m-1}{2}n} 2^{mn} b^{n(mn+1)}$$

(5.9)

*and $dd_0^{(m)}(n,0,b) = 0$ for all other $n$.*

For $a = 0$ the determinants $d_0^{(m)}(n,0,b)$ and $d_1^{(m)}(n,0,b)$ can also be computed with the **continued fractions method** by **Gessel** and **Xin** ([9],[12]). Suppose that $D(x,y) = \sum_{i,j=0}^{\infty} d_{i,j} x^i y^j$ and let $[D(x,y)]_n = \det(d_{i,j})_{i,j=0}^{n-1}$. If $u(x)$ is a formal power series with $u(0) = 1$ then $[u(x)D(x,y)]_n = [D(x,y)]_n$ because the transformed determinant is obtained from the original determinant by elementary row operations.

Then the Hankel determinants $H_n(A) = \det(a(i+j))_{i,j=0}^{n-1}$ and $H_n^1(A) = \det(a(i+j+1))_{i,j=0}^{n-1}$ of the formal power series $A(z) = \sum_{n\geq 0} a(n) z^n$ are given by

$$H_n(A) = \left[\frac{xA(x) - yA(y)}{x - y}\right]_n \text{ and } H_n^1(A) = \left[\frac{A(x) - A(y)}{x - y}\right]_n.$$

Let now $A_0(x) = f_m(x,0,b) = \dfrac{1}{1 - bx^m A_0(x)}$.

Then
$$H_n(A_0) = \left[\frac{1}{x-y}\left(\frac{x}{1-bx^m A_0(x)} - \frac{y}{1-by^m A_0(y)}\right)\right]_n = \left[\frac{x - y + bxy\left(xx^{m-2}A_0(x) - yy^{m-2}A_0(y)\right)}{(1-bx^m A_0(x))(1-by^m A_0(y))(x-y)}\right]_n$$

$$= \left[1 + xy\frac{xA_1(x) - yA_1(y)}{x-y}\right]_n = H_{n-1}(A_1)$$

with $A_1(x) = bx^{m-2} A_0(x)$ which satisfies $A_1(x) = \dfrac{bx^{m-2}}{1 - x^2 A_1(x)}$.

Let $A_2(x) = \dfrac{1}{b} A_1(x) = \dfrac{x^{m-2}}{1 - bx^2 A_2(x)}$. Then $H_n(A_1) = b^n H_n(A_2)$.

Now

$$H_n(A_2) = \left[ \frac{1}{x-y} \left( \frac{x^{m-1}}{1-bx^2 A_2(x)} - \frac{y^{m-1}}{1-by^2 A_2(y)} \right) \right]_n =$$

$$\left[ \frac{1}{x-y} \left( x^{m-1} - y^{m-1} + (xy)^{m-1} \right) \left( \frac{xA_3(x) - yA_3(y)}{x-y} \right) \right]_n$$

$$= (-1)^{\binom{m-1}{2}} H_{n-m+1}(A_3)$$

with $A_3(x) = \dfrac{b}{x^{m-2}} A_2(x) = \dfrac{b}{1-x^2 A_2(x)}$.

Let $A_4(x) = \dfrac{A_3(x)}{b} = \dfrac{1}{1-bx^m A_4(x)} = A_0(x)$.

Thus $H_n(A_3) = b^n H_n(A_0)$.

Therefore we get $H_n(A_0) = b^{n-1} H_{n-1}(A_2) = (-1)^{\binom{m-1}{2}} b^{n-1} H_{n-m}(A_3) = (-1)^{\binom{m-1}{2}} b^{2n-m-1} H_{n-m}(A_0)$.

It is easily verified that $H_1(A_0) = d_0^{(m)}(1,0,b) = 1$, $H_m(A_0) = d_0^{(m)}(m,0,b) = (-1)^{\binom{m-1}{2}} b^{m-1}$ and $H_n(A_0) d_0^{(m)}(n,0,b) = 0$ for $1 < n < m$.

With the same method we get

**Theorem 5.6**

*For $m \geq 1$*

$$d_1^{(m)}(mn,0,b) = (-1)^{\binom{m}{2}n} b^{mn^2}. \tag{5.10}$$

$d_1^{(m)}(n,0,b) = 0$ *else*.

**Proof**

$$H_n^1(A_0) = \left[ \frac{1}{x-y} \left( \frac{1}{1-bx^m A_0(x)} - \frac{1}{1-by^m A_0(y)} \right) \right]_n = \left[ \frac{bx^m A_0(x) - by^m A_0(y)}{(x-y)} \right]_n$$

$$= \left[ \frac{xA_1(x) - yA_1(y)}{x-y} \right]_n = H_n(A_1)$$

with $A_1(x) = bx^{m-1} A_0(x)$ or $A_1(x) = \dfrac{bx^{m-1}}{1-xA_1(x)}$.

Let $A_2(x) = \dfrac{x^{m-1}}{1-bxA_2(x)}$. Then $H_n(A_1) = b^n H_n(A_2)$.

$$H_n(A_2) = \left[\frac{1}{x-y}\left(\frac{x^m}{1-bxA_2(x)} - \frac{y^m}{1-byA_2(y)}\right)\right]_n = \left[\frac{1}{x-y}\left(x^m - y^m + (xy)^m(A_3(x) - A_3(y))\right)\right]_n$$

$$= (-1)^{\binom{m}{2}} H^1_{n-m}(A_3)$$

with $A_3(x) = \dfrac{b}{x^{m-1}} A_2(x)$ or $A_3(x) = \dfrac{b}{1 - x^m A_3(x)}$.

Then $A_3 = bA_0$ and $H^1_n(A_3) = b^n H^1_n(A_0)$.

This implies $H^1_n(A_0) = (-1)^{\binom{m}{2}} b^{2n-m} H^1_{n-m}(A_0)$.

Since $H^1_n(A_0) = 0$ for $n < m$ and $H^1_m(A_0) = (-1)^{\binom{m}{2}} b^m$ we get by induction Theorem 5.6.

In the general case we use the **Lindström–Gessel-Viennot theorem** (cf. e.g. [10], Theorem 6, [3], [12],[14]). It implies that

$$d_k^{(m)}(n,a,b) = \sum_\sigma \operatorname{sgn}\sigma \cdot W(\sigma) \tag{5.11}$$

where $\sigma$ runs through all permutations of $\{0,1,\cdots,n-1\}$ and $W(\sigma)$ is the weight of all systems $S_k(n)$ of $n$ non-intersecting non-negative lattice paths from $A_i \to E_{\sigma(i)}$ for all $i \in \{0,1,\cdots,n-1\}$.

Here the initial points are $A_i = (-i,0)$ and for given $k$ the endpoints are $E_i = (i+k,0)$. From the context it will always be clear which $k$ is considered. So we do not indicate $k$ explicitly.

A set of $n$ non-overlapping lattice paths from $\{A_0,\cdots,A_{n-1}\} \to \{E_0,\cdots,E_{n-1}\}$ will be called $n-$**set**.

**6. The case $m=3$.**

Let us first consider the situation for $m=3$.

Let $S$ be an $n-$set with permutation $\sigma$. Let $w_0(S)$ be the weight of $S$ and let the signed weight $w(S)$ be $w(S) = \operatorname{sgn}(\sigma) w_0(S)$.

a) The Hankel determinants $d_0^{(3)}(n,a,b)$.

First some remarks on notation. For $k=0$ we denote by $Id_0$ the trivial path $A_0 = (0,0) \to E_0 = (0,0)$. Each other path consists of a sequence of steps $U, H, D$ from left to right. If the path begins in $A_i$ the first step will get the index $i$, and if it ends in $E_j$ then the last step will get the index $j$. For example for $k = 0$ and $m = 3$ the path $A_3 = (-3,0) \to (-2,1) \to (-1,1) \to (0,1) \to (2,0) = E_2$ will be denoted by $U_3 H^2 D_2$.

For $k = 2$ the path $A_2 = (-2,0) \to (-1,1) \to (0,2) \to (2,1) \to (3,1) \to (5,0) = E_3$ is denoted by $U_2^2 D H D_3$.

**Theorem 6.1**

$$d_0^{(3)}(3n,a,b) = (-1)^n b^{n(3n-1)},$$
$$d_0^{(3)}(3n+1,a,b) = (-1)^n b^{n(3n+1)}, \qquad (6.1)$$
$$d_0^{(3)}(3n+2,a,b) = 0.$$

*Furthermore we have* $D_0^{(3)}(n,a,b,t) = d_0^{(3)}(n,a,b).$

**Proof**

For $n = 1$ we have the trivial $1-$ set $\{Id_0\}$.

For $n = 2$ there is no $2-$ set because no path which starts at $A_1$ can end in $E_1$.

For $n \geq 3$ each $n-$ set must contain the $3-$ set $S_3 = \{Id_0, U_2 D_1, U_1 D_2\}$. The signed weight is $w(S_3) = -b^2$.

If $i \geq 3$ all paths starting from $A_i$ must begin with two up-steps and those ending in $E_i$ end with two down-steps. This implies that on height $2$ we have the same situation as on height $0$ but with points $A_{i-3}^* = E_{i-3}^*$ in place of $A_0 = E_0$. For an $n-$ set $S$ the weight of the part of all paths which lie between height $y = 0$ and height $y = 2$ is $-b^2 b^{2(n-3)} = -b^{2n-4}$.

Therefore we get

$$d_0^{(3)}(n,a,b) = -b^{2n-4} d_0^{(3)}(n-3,a,b) \qquad (6.2)$$

with initial values $d_0^{(3)}(1,a,b) = 1$ and $d_0^{(3)}(2,a,b) = 0$. Since $d_0^{(3)}(3,a,b) = -b^2 = -b^{2 \cdot 3 - 4}$ we can set $d_0^{(3)}(0,a,b) = 1$.

Thus by induction we get $d_0^{(3)}(3n,a,b) = -b^{6n-4}(-1)^{n-1} b^{(n-1)(3n-4)} = (-1)^n b^{n(3n-1)}$ and
$d_0^{(3)}(3n+1,a,b) = -b^{6n-2}(-1)^{n-1} b^{(n-1)(3n-2)} = (-1)^n b^{n(3n+1)}.$

The same argument holds for $D_0^{(3)}(n,a,b,t)$.

b) In an analogous way we can compute $d_1^{(3)}(n,a,b)$.

**Theorem 6.2**

$$d_1^{(3)}(3n,a,b) = (-1)^n b^{3n^2},$$
$$d_1^{(3)}(3n+1,a,b) = (-1)^n (n+1)ab^{3n^2+2n}, \qquad (6.3)$$
$$d_1^{(3)}(3n-1,a,b) = (-1)^{n-1} nab^{3n^2-2n}.$$

**Proof**

We get the recursion

$$d_1^{(3)}(n,a,b) = ab^{n-1} d_0^{(3)}(n-1,a,b) - b^{2n-3} d_1^{(3)}(n-3,a,b). \qquad (6.4)$$

If there is a path $A_0 = (0,0) \to (1,0) = E_0$ in the $n-$set $S_{n,0}$ then on height 1 we have the points $A_{n-2}^* = (-n+2,1), \cdots, A_0^* = E_0^* = (0,1), E_1^* = (1,1), \cdots E_{n-2}^* = (n-2,1)$. On these points we can construct an arbitrary $(n-1)-$set for $d_0^{(3)}(n-1,a,b)$. The signed weight of the paths till height 1 is $ab^{n-1}$.

If there is no path $(0,0) \to (1,0)$, then there must be paths 
$A_2 = (-2,0) \to (-1,1) \to (1,0) = E_0, A_1 = (-1,0) \to (0,1) \to (2,0) = E_1, A_0 = (0,0) \to (1,1) \to (3,0) = E_2$. This gives the $2-$set $\{U_0 D_2, U_1 D_1, U_2 D_0\}$.

On height 2 we have the points $A_{n-3}^{**} = (-n+3,2), \cdots, A_0^{**} = (-1,2), E_0^{**} = (0,2), \cdots, E_{n-3}^{**} = (n-4,2)$. On these points we can construct an arbitrary $(n-3)-$set for $d_1^{(3)}(n-3,a,b)$. The signed weight of the paths till height 2 is $-b^{2n-3}$. This proves (6.4).

Now we can by induction prove the explicit formulae:

$$d_1^{(3)}(3n,a,b) = ab^{3n-1} d_0^{(3)}(3n-1,a,b) - b^{6n-3} d_1^{(3)}(3n-3,a,b) = -b^{6n-3}(-1)^{n-1} b^{3(n-1)^2} = (-1)^n b^{3n^2}.$$

$$d_1^{(3)}(3n+1,a,b) = ab^{3n} d_0^{(3)}(3n,a,b) - b^{6n-1} d_1^{(3)}(3n-2,a,b)$$
$$= (-1)^n ab^{3n} b^{n(3n-1)} - b^{6n-1}(-1)^{n-1} nab^{3n^2-4n+1} = (-1)^n (n+1)ab^{3n^2+2n}$$

$$d_1^{(3)}(3n+2,a,b) = ab^{3n+1} d_0^{(3)}(3n+1,a,b) - b^{6n+1} d_1^{(3)}(3n-1,a,b)$$
$$= (-1)^n ab^{3n+1} b^{n(3n+1)} - ab^{6n+1}(-1)^{n-1} nb^{3(n-1)^2+4(n-1)+1} = (-1)^n(n+1)ab^{3n^2+4n+1}.$$

To compute $D_1^{(3)}(n,a,b,t)$ we note that the only difference consists in the fact that on height 0 the weight of $H$ is $a+t$ instead of $a$. Therefore we get

$$D_1^{(3)}(n,a,b) = (a+t)b^{n-1} d_0^{(3)}(n-1,a,b) - b^{2n-3} d_1^{(3)}(n-3,a,b) = d_1^{(3)}(n,a,b) + tb^{n-1} d_0^{(3)}(n-1,a,b).$$

This gives

**Theorem 6.3**

$$D_1^{(3)}(3n,a,b,t) = (-1)^n b^{3n^2},$$
$$D_1^{(3)}(3n+1,a,b,t) = (-1)^n ((n+1)a+t)b^{3n^2+2n}, \quad (6.5)$$
$$D_1^{(3)}(3n-1,a,b,t) = (-1)^{n-1}(na+t)b^{3n^2-2n}.$$

c) The determinants $d_2^{(3)}(n,a,b)$.

**Theorem 6.4**

$$d_2^{(3)}(3n,a,b) = (-1)^{n+1} b^{3n^2+n-1}\left(a^3 \sum_{i=0}^{n} i^2 - (n+1)b\right),$$
$$d_2^{(3)}(3n+1,a,b) = (-1)^n (n+1)^2 a^2 b^{3n^2+3n}, \quad (6.6)$$
$$d_2^{(3)}(3n+2,a,b) = (-1)^n b^{3n^2+5n+1}\left(a^3 \sum_{i=0}^{n+1} i^2 - (n+1)b\right).$$

**Proof**

It would be nice if we could prove this too with the help of the condensation formula

$$d_2(n)d_0(n) = d_0(n+1)d_2(n-1) + d_1(n)^2.$$

But due to the fact that $d_0^{(3)}(3n-1,a,b) = 0$ the condensation formula only gives

$$d_2^{(3)}(3n+1,a,b) = (-1)^n (n+1)^2 a^2 b^{3n^2+3n} \quad (6.7)$$

and

$$d_2^{(3)}(3n,a,b) = b^{2n} d_2^{(3)}(3n-1,a,b) + (-1)^n b^{3n^2+n}. \quad (6.8)$$

To prove Theorem 6.4 we need the extra information

$$d_2^{(3)}(3n+2,a,b) = ab^{2n+1} d_2^{(3)}(3n+1,a,b) - b^{4n+2} d_2^{(3)}(3n,a,b). \quad (6.9)$$

It is clear that each $3n-$ set is contained in the region beneath the line segments
$(-3n+1,0) \to (-n+1,2n)$ and $(-n+1,2n) \to (3n+1,0)$,

each $(3n+1)-$ set beneath $(-3n,0) \to (-n,2n)$, $(-n,2n) \to (-n+2,2n)$ and $(-n+2,2n) \to (3n+2,0)$,

and each $(3n+2)-$ set beneath $(-3n-1,0) \to (-n,2n+1)$, $(-n,2n+1) \to (-n+1,2n+1)$ and
$(-n+1,2n+1) \to (3n+3,0)$.

The heights are exact with the exception of the $3n$ – sets $S^*_{3n}$ defined by

$$S^*_{3n} = \{S^*_{3(n-1)}, U^{2n-1}_{3n-3}D^{2n-1}_{3n-2}, U^{2n-1}_{3n-2}D^{2n-1}_{3n-3}, U^{2n-1}_{3n-1}DUD^{2n-1}_{3n-1}\},$$

which are contained in the smaller region bounded by

$(-3n+1, 0) \to (-n, 2n-1)$, $(-n+2, 2n-2) \to (-n+3, 2n-1)$, and $(-n+3, 2n-1) \to (3n+1, 0)$.

We call these $3n$ – sets exceptional and the other $3n$ – sets normal.

It is easily seen that no $(3n+1)$ – set can contain a normal $3n$ – set and that each exceptional $S^*_{3n}$ is contained in precisely one $(3n+1)$ – set $\{S^*_{3n}, U^{2n}_{3n}H^2D^{2n}_{3n}\}$ and in precisely one $(3n+2)$ – set $\{S^*_{3n}, U^{2n+1}_{3n}D^{2n+1}_{3n+1}, U^{2n+1}_{3n+1}D^{2n+1}_{3n}\}$.

Each normal $3n$ – set is also contained in precisely one $(3n+2)$ – set.

Each $(3n+1)$ – set can also uniquely be enlarged to a $(3n+2)$ – set by adjoining the path $U^{2n+1}_{3n+1}HD^{2n+1}_{3n+1}$.

This implies that the signed weight $d^{(3)}_2(3n+2, a, b)$ of all $(3n+2)$ – sets satisfies

$$d^{(3)}_2(3n+2, a, b) = ab^{2n+1}d^{(3)}_2(3n+1, a, b) - b^{4n+2}d^{(3)}_2(3n, a, b).$$

It is now routine matter to verify (6.6).

To compute the Hankel determinants $D^{(3)}_2(n, a, b, t)$ we get from the condensation formula

$$D^{(3)}_2(3n+1, a, b, t) = (-1)^n((n+1)a+t)^2 b^{3n^2+3n} \tag{6.10}$$

and

$$D^{(3)}_2(3n, a, b, t) = b^{2n}D^{(3)}_2(3n-1, a, b, t) + (-1)^n b^{3n^2+n}. \tag{6.11}$$

We need the extra information

$$D^{(3)}_2(3n+2, a, b, t) = ab^{2n+1}D^{(3)}_2(3n+1, a, b, t) - b^{4n+2}D^{(3)}_2(3n, a, b, t). \tag{6.12}$$

This follows in the same way as above.

Now we can deduce

**Theorem 6.5**

$$D_2^{(3)}(3n,a,b,t) = (-1)^n b^{3n^2+n-1}\left(a^3\sum_{i=0}^{n} i^2 - (n+1)b + nat(t+(n+1)a)\right),$$

$$D_2^{(3)}(3n+1,a,b,t) = (-1)^n (t+(n+1)a)^2 b^{3n^2+3n}, \qquad (6.13)$$

$$D_2^{(3)}(3n-1,a,b,t) = (-1)^{n-1} b^{3n^2-n-1}\left(a^3\sum_{i=0}^{n} i^2 - nb + nat(t+(n+1)a)\right).$$

**Another formulation**

From (6.5) we see that the sequence $\left(\dfrac{C(n+1,0;3,a,b,t)}{a}\right)_{n\geq 0}$ satisfies (2.1). Therefore there exist uniquely determined numbers $s(n)$ and $t(n)$ which satisfy (2.5).

**Theorem 6.6**

*The uniquely determined numbers $s(n)$ and $t(n)$ for the sequence $\left(\dfrac{C(n+1,0;3,a,b,t)}{a}\right)_{n\geq 0}$ are*

$$s(3n) = (n+1)a + t, \; s(3n+1) = a, \; s(3n+2) = -(n+1)a - t \qquad (6.14)$$

*and*

$$t(3n) = \frac{b}{(n+1)a+t}, \; t(3n+1) = -\frac{b}{(n+1)a+t}, \; t(3n+2) = -\big((n+1)a+t\big)\big((n+2)a+t\big). \qquad (6.15)$$

**Proof**

The sequence $(t(n))_{n\geq 0}$ can be computed from (2.6). This gives

$$\prod_{i=1}^{n-1}\prod_{j=0}^{i-1} t(j) = D_1^{(3)}(n,a,b,t) \text{ or } t(n) = \frac{D_1^{(3)}(n,a,b,t) D_1^{(3)}(n+2,a,b,t)}{D_1^{(3)}(n+1,a,b,t)^2},$$

which implies (6.15).

From (2.7) we deduce that $(-1)^n p(n,0) = \dfrac{D_2^{(3)}(n,a,b,t)}{D_1^{(3)}(n,a,b,t)}.$

This implies that

$$p(3n,0) = (-1)^{n+1} b^{n-1} \left( a^3 \sum_{i=0}^{n} i^2 - (n+1)b + nat(t + (n+1)a) \right),$$

$$p(3n+1,0) = (-1)^{n+1} (t + (n+1)a) b^n,$$

$$p(3n+2,0) = \frac{(-1)^n b^n \left( a^3 \sum_{i=0}^{n+1} i^2 - (n+1)b + (n+1)at(t + (n+2)a) \right)}{(n+1)a + t}.$$

By (2.3) we have $p(n,x) = (x - s(n-1))p(n-1,x) - t(n-2)p(n-2,x)$. Therefore
$p(n,0) = -s(n-1)p(n-1,0) - t(n-2)p(n-2,0)$.

This implies that

$$s(n) = -\frac{p(n+1,0) + t(n-1)p(n-1,0)}{p(n,0)}$$

and we get (6.14).

I could not find a direct proof of (6.14) and (6.15).

For the case of unrestricted paths I have only some conjectures based on computer experiments.

**Conjecture 6.7**

Let $\dfrac{1}{\sqrt{(1-az)^2 - 4bz^3}} = \sum_{n \geq 0} g(n,3,a,b) z^n$. The uniquely determined numbers $s(n)$ and $t(n)$ for the sequence $\left( \dfrac{g(n+1,3,a,b)}{a} \right)_{n \geq 0}$ are

$s(0) = a$, $s(3n) = \dfrac{2n+1}{2} a$ for $n > 0$, $s(3n+1) = a$, $s(3n+2) = -\dfrac{2n+1}{2} a$

and

$$t(3n) = \frac{2b}{(2n+1)a}, \quad t(3n+1) = -\frac{2b}{(2n+1)a}, \quad t(3n+2) = -\frac{a^2}{4}(2n+1)(2n+3).$$

This gives

**Conjecture 6.8**

$$dd_0^{(3)}(3n,a,b) = (-1)^n 2^{3n-1} b^{n(3n-1)},$$
$$dd_0^{(3)}(3n+1,a,b) = (-1)^n 2^{3n} b^{n(3n+1)}, \qquad (6.16)$$
$$dd_0^{(3)}(3n+2,a,b) = 0.$$

$$dd_1^{(3)}(3n,a,b) = (-1)^n 2^{3n} b^{3n^2},$$
$$dd_1^{(3)}(3n+1,a,b) = (-1)^n (2n+1) 2^{3n} ab^{3n^2+2n} \qquad (6.17)$$
$$dd_1^{(3)}(3n-1,a,b) = (-1)^{n-1} (2n-1) 2^{3n-2} ab^{3n^2-2n}$$

$$dd_2^{(3)}(3n+1,a,b) = (-1)^n (2n+1)^2 2^{3n} a^2 b^{3n(n+1)},$$
$$dd_2^{(3)}(3n,2,a,b) = (-1)^n \left( (2n+1) 2^{3n} b^{n(3n+1)} - \binom{2n+1}{3} 2^{3n-1} a^3 b^{3n^2+n-1} \right), \qquad (6.18)$$
$$dd_2^{(3)}(3n-1,2,a,b) = (-1)^n \left( (2n-1) 2^{3n-1} b^{n(3n-1)} - \binom{2n+1}{3} 2^{3n-2} a^3 b^{3n^2-n-1} \right).$$

**Remark 6.10**

(6.17) follows immediately from Conjecture 6.7.

To deduce (6.18) let $u(n) = s(n-1)u(n-1) - t(n-2)u(n-2)$ with initial values $u(0) = 1$ and $u(1) = s(0)$.

It is easily verified that

$$u(3n+1) = (2n+1)ab^n,$$
$$u(3n) = (2n+1)b^n - \frac{a^3}{2} b^{n-1} \binom{2n+1}{3},$$
$$u(3n-1) = -\frac{2}{a} b^n + \frac{a^2 b^{n-1}}{3} \binom{2n+1}{2}.$$

Using (2.7) this implies (6.18).

The condensation formula implies (6.16).

**An example**

Consider the sequence $(c(n,3,1,1))_{n \geq 0} = (1,1,1,2,4,7,13,26,\cdots)$ (OEIS A023431).

The Hankel determinants are

$$\left(d_0^{(3)}(3n,1,1), d_0^{(3)}(3n+1,1,1), d_0^{(3)}(3n+2,1,1)\right)_{n\geq 0} = \left((-1)^n, (-1)^n, 0\right)_{n\geq 0},$$

$$\left(d_1^{(3)}(3n,1,1), d_1^{(3)}(3n+1,1,1), d_1^{(3)}(3n+2,1,1)\right)_{n\geq 0} = \left((-1)^n, (-1)^n(n+1), (-1)^n(n+1)\right)_{n\geq 0},$$

$$\left(d_2^{(3)}(3n,1,1), d_2^{(3)}(3n+1,1,1), d_2^{(3)}(3n+2,1,1)\right)_{n\geq 0} = \left((-1)^{n+1}\binom{n+2}{2}\frac{2n-3}{3}, (-1)^n(n+1)^2, (-1)^n\binom{n+1}{2}\frac{2n+7}{3}\right).$$

For the sequence $(g(n,3,1,1)) = (1,1,1,3,7,13,27,61,\cdots)$ (OEIS A01084710) we get

$$\left(dd_0^{(3)}(3n,1,1), dd_0^{(3)}(3n+1,1,1), dd_0^{(3)}(3n+2,1,1)\right)_{n\geq 0} = \left((-1)^n 2^{3n-1}, (-1)^n 2^{3n}, 0\right)_{n\geq 0},$$

$$\left(dd_1^{(3)}(3n,1,1), dd_1^{(3)}(3n+1,1,1), dd_1^{(3)}(3n+2,1,1)\right)_{n\geq 0} = \left((-1)^n 2^{3n}, (-1)^n(2n+1)2^{3n}, (-1)^n(2n+1)2^{3n+1}\right)_{n\geq 0},$$

$$\left(dd_2^{(3)}(3n,1,1), dd_2^{(3)}(3n+1,1,1), dd_2^{(3)}(3n+2,1,1)\right)_{n\geq 0}$$

$$= \left((-1)^n 2^{3n-1}\left(4n+2-\binom{2n+1}{3}\right), (-1)^n 2^{3n}(2n+1)^2, (-1)^n 2^{3n+1}\left(4n+2-\binom{2n+3}{3}\right)\right).$$

## 7. The general case

**Theorem 7.1**

*Let $m \geq 2$. Then*

$$d_0^{(m)}(mn,a,b) = (-1)^{\binom{m-1}{2}n} b^{n(mn-1)}$$

$$d_0^{(m)}(mn+1,a,b) = (-1)^{\binom{m-1}{2}n} b^{n(mn+1)}$$

(7.1)

*and $d_0^{(m)}(n,a,b) = 0$ for all other $n$.*

**Proof**

The proof is almost the same as for $m = 3$.

For $n = 1$ we have only the trivial $1-$ set $\{Id_0\}$.

For $1 < n < m$ there is no $n-$ set since no path originating in $A_1$ is contained in $\{1, 2, \cdots, n-1\}$.

For $n = m$ there is one $n-$ set $S_m = \{Id_0, U_1D_{m-1}, U_2D_{m-2}, \cdots, U_{m-1}D_1\}$ with $w(S_m) = (-1)^{\binom{m-1}{2}} b^{m-1}$.

For $n \geq m$ each $n-$ set must contain $S_m$ together with the paths $U_i^2$ and $D_i^2$ for $m \leq i \leq n-1$. On height 2 we have the same situation as in the beginning with $n$ replaced by $n-m$. The weight of the paths till height 2 is $b^{m-1}b^{2(n-m)} = b^{2n-m-1}$.

Therefore we get by induction the above results.

Now we want to compute $d_1^{(m)}(n,a,b)$.

For $n=1$ we have the $1-$ set $S_1 = \{H_0\} = \{A_0 \to E_0\}$ with signed weight $w(S_1) = d_1^{(m)}(0,a,b) = a$.

First we want to construct all $n-$ sets $S_n$ which contain $S_1$.

In this case each of the points $A_i$, $i \geq 1$, must be the begin of an up-step $U$ and each of the points $E_i$, $i \geq 1$, must be the end a down-step. The smallest such $n-$ set $S_{m-1} = \{S_1, U_1 D_{m-2}, U_2 D_{m-2}, \cdots, U_{m-2} D_1\}$ with weight $w(S_{m-1}) = -(-1)^{\binom{m}{2}} ab^{m-2}$ occurs for $n = m-1$.

For $i \geq m-1$ all paths starting from $A_i$ must begin with two up-steps $U^2$ and all paths with endpoint $E_i$ must end with at least $2$ down-steps $D^2$.

The 2 paths $U_{m-1}^2 D_m^2 : A_{m-1} \to E_m : (-m+1, 0) \to (-m+3, 2) \to (m+1, 0)$ and $U_m^2 D_{m-1}^2 : A_m \to E_{m-1} : (-m, 0) \to (-m+2, 2) \to (m, 0)$ extending $S_{m-1}$ give the $(m+1)-$ set $S_{m+1}$ with weight $w(S_{m+1}) = (-1)^{\binom{m}{2}} ab^{m+2}$.

Let now $u(n,m,a,b) = w(S_n)$ if there is an $n-$ set $S_n$ which contains $S_1$ and $u_n = 0$ else.

Thus for $n > m+1$ on height $2$ we have the same situation as after the first step: The neighbouring points $A_0^* = (-m+2, 2)$ and $E_0^* = (-m+3, 2)$ are occupied and we look for non-overlapping paths between the points $A_i^* = (-m+2-i, 2)$ and $E_j^* = (-m+3+j, 2)$ for $1 \leq i, j \leq n-m-1$. This gives

$$u(n,m,a,b) = (-1)^{\binom{m}{2}} u(n-m, m, a, b) b^{2n-m}$$ with initial values $u(0,m,a,b) = 0, u(1,m,a,b) = a,$

$u(n,m,a,b) = 0$ for $1 < n < m-1$, $u(m-1, m, a, b) = -(-1)^{\binom{m}{2}} ab^{m-2}$. (Observe that

$u(m+1, m, a, b) = (-1)^{\binom{m}{2}} b^{2(m+1)-m} a.$)

From this we get by induction

$$u(mn-1, m, a, b) = -(-1)^{n\binom{m}{2}} ab^{mn^2 - 2n}$$
$$u(mn+1, m, a, b) = (-1)^{n\binom{m}{2}} ab^{mn^2 + 2n}$$
(7.2)

and $u(n,m,a,b) = 0$ else.

There remains the case where there is no path $A_0 \to E_0$.

Here we get the smallest $n$-set $T_n$ for $n = m$.

Let $v(n,m,a,b) = w(T_n)$ if there is such an $n$-set $T_n$ and $v_n = 0$ else.

We have $T_m = \{U_0 D_{m-1}, U_1 D_{m-2}, \cdots, U_{m-1} D_0\}$ with signed weight $v(m,m,a,b) = w(T_m) = (-1)^{\binom{m}{2}} b^m$.

Let now $n > m$. For $i = m$ we get on height $2$ the point $(-m+2, 2)$ as image of $A_m$ and the point $(-m+3, 2)$ as pre-image of $E_m$. Therefore we have the same situation as in the definition of $d_1^{(m)}(n-m, a, b)$ and get

$$v(n,m,a,b) = (-1)^{\binom{m}{2}} b^{2n-m} d_1^{(m)}(n-m, a, b).$$

This gives with induction

**Theorem 7.2**

$$d_1^{(m)}(mn, a, b) = (-1)^{\binom{m}{2} n} b^{mn^2}.$$

$$d_1^{(m)}(mn+1, a, b) = (-1)^{\binom{m}{2} n} (n+1) a b^{mn^2 + 2n}$$

$$d_1^{(m)}(mn-1, a, b) = -(-1)^{\binom{m}{2} n} n a b^{mn^2 - 2n}$$

$d_1^{(m)}(n, a, b) = 0$ else for $m \geq 3$.

**Proof**

This result holds for $n = 0$. If it holds for $n - 1$ then we get

$$d_1^{(m)}(m(n-1)+1, a, b) = (-1)^{\binom{m}{2}(n-1)} n a b^{m(n-1)^2 + 2(n-1)}$$

and therefore

$$d_1^{(m)}(mn+1, a, b) = (-1)^{n\binom{m}{2}} a b^{mn^2 + 2n} + (-1)^{\binom{m}{2}} b^{2(mn+1)-m} (-1)^{\binom{m}{2}(n-1)} n a b^{m(n-1)^2 + 2(n-1)}$$
$$= (-1)^{n\binom{m}{2}} (n+1) a b^{mn^2 + 2n}.$$

In the same way we get

$$d_1^{(m)}(mn-1,a,b) = (-1)^{n\binom{m}{2}} ab^{mn^2-2n} + (-1)^{\binom{m}{2}} b^{2(mn-1)-m} (-1)^{\binom{m}{2}(n-1)} nab^{m(n-1)^2-2(n-1)}$$

$$= (-1)^{n\binom{m}{2}} (n+1) ab^{mn^2-2n}.$$

**Remark 7.3**

Let $m \geq 2$ and

$$\tilde{p}_n(m,z) := \det \begin{pmatrix} a(0) & a(1) & \cdots & a(n-1) & 1 \\ a(1) & a(2) & \cdots & a(n) & z \\ a(2) & a(3) & \cdots & a(n+1) & z^2 \\ \vdots & & & & \vdots \\ a(n) & a(n+1) & \cdots & a(2n-1) & z^n \end{pmatrix} \quad (7.3)$$

with $a(n) = c(n,m,a,b)$. Since in general the coefficient of $z^n$ can vanish we have no direct analogue of $p_n(z)$. But if we define

$$p_n(m,z) = \frac{1}{b^{\frac{n(n-1)}{m}}} \tilde{p}_n(m,z) \quad (7.4)$$

then we have $p_n(2,z) = p_n(z)$.

We conjecture that

$$\tilde{p}_{mn}(m,z) = (-1)^{n\binom{m-1}{2}} b^{n(mn-1)} Fib_{2n+1}\left(\sqrt{z^{m-2}}(z-a), -b\right),$$

$$\tilde{p}_{mn+1}(m,z) = (-1)^{n\binom{m-1}{2}} b^{n(mn+1)} z^{1-\frac{m}{2}} Fib_{2n+2}\left(\sqrt{z^{m-2}}(z-a), -b\right), \quad (7.5)$$

$$\tilde{p}_{mn+m-1}(m,z) = (-1)^{\binom{m-1}{2}} b^{(m-2)(2n+1)} \tilde{p}_{mn+1}(m,z)$$

and $\tilde{p}_n(m,z) = 0$ else.

For $m = 2$ this coincides with (2.14).

Note that $Fib_n((z-a)\sqrt{z^{m-2}}, -b) = z^{\frac{n(2-m)}{2}} Fib_n((z-a), -bz^{2-m})$.

This conjecture is equivalent with

$$p_n(m,z) = (-1)^{\binom{m-1}{2}} \left(z^{m-2}(z-a)^2 - 2b\right) p_{n-m}(m,z) - b^2 p_{n-2m}(m,z) \quad (7.6)$$

for $n \geq 2m$ with initial values

$$p_0(m,z) = 1, \; p_1(m,z) = z - a, \; p_{m-1}(m,z) = (-1)^{\binom{m-1}{2}}(z-a)b^{\frac{m-2}{m}},$$

$$p_m(m,z) = (-1)^{\binom{m-1}{2}}\left(z^{m-2}(z-a)^2 - b\right), \; p_{m+1}(m,z) = (-1)^{\binom{m-1}{2}}(z-a)\left(z^{m-2}(z-a)^2 - 2b\right),$$

$$p_{2m-1}(m,z) = b^{\frac{m-2}{m}}(z-a)\left(z^{m-2}(z-a)^2 - 2b\right)$$

and $p_n(m,z) = 0$ for the other values of $n$ with $1 < n < 2m$.

For the determinants $D_1^{(m)}(n,a,b,t)$ we get with the same reasoning

$$D_1^{(m)}(n,a,b,t) = u(n,m,a+t,b) + v(n,m,a,b) = u(n,m,a,b) + u(n,m,t,b) + v(n,m,a,b)$$
$$= d_1^{(m)}(n,a,b) + tu(n,m,1,b)$$

Taking into account (7.2) we get

**Theorem 7.3**

$$D_1^{(m)}(mn,a,b,t) = (-1)^{\binom{m}{2}n} b^{mn^2}.$$

$$D_1^{(m)}(mn+1,a,b,t) = (-1)^{\binom{m}{2}n}(t+(n+1)a)b^{mn^2+2n}$$

$$D_1^{(m)}(mn-1,a,b,t) = -(-1)^{\binom{m}{2}n}(t+na)b^{mn^2-2n}$$

$D_1^{(m)}(n,a,b,t) = 0$ else for $m \geq 3$.

**Theorem 7.4**

*For $m \geq 4$ the following identities hold:*

$$d_2^{(m)}(mn+1,a,b) = (-1)^{n\binom{m-1}{2}}(n+1)^2 a^2 b^{mn^2+3n},$$
$$d_2^{(m)}(mn-2,a,b) = -(-1)^{n\binom{m-1}{2}} n^2 a^2 b^{mn^2-3n}. \tag{7.7}$$

**Proof**

The condensation formula $d_2(n)d_0(n) = d_0(n+1)d_2(n-1) + d_1(n)^2$ gives

$$d_2^{(m)}(mn+1,a,b) = (-1)^{n\binom{m-1}{2}}(n+1)^2 a^2 b^{mn^2+3n},$$

$$d_2^{(m)}(mn-2,a,b) = -(-1)^{n\binom{m-1}{2}} n^2 a^2 b^{mn^2-3n},$$

and

$$d_2^{(m)}(mn,a,b) = d_2^{(m)}(mn-1,a,b)b^{2n} + (-1)^{n\binom{m-1}{2}} b^{mn^2+n}. \tag{7.8}$$

But we have

**Conjecture 7.5**

$$d_2^{(m)}(mn,a,b) = (-1)^{n\binom{m-1}{2}}(n+1)b^{mn^2+n},$$

$$d_2^{(m)}(mn-1,a,b) = (-1)^{n\binom{m-1}{2}} n b^{mn^2-n}.$$

This follows from the conjectured identity

$$d_2^{(m)}(mn,a,b) = (-1)^{\binom{m-1}{2}} b^{2nm-m+1} d_2^{(m)}(mn-m,a,b) + (-1)^{n\binom{m-1}{2}} b^{mn^2+n} \tag{7.9}$$

or from the conjecture that for $m \geq 4$ no $mn-$ set can contain a horizontal step.

I suppose that there must be a simple proof of this fact, but unfortunately I could not find it.

In the same way we obtain

**Conjecture 7.6**

*For $m \geq 4$ the following identities hold:*

$$D_2^{(m)}(mn,a,b,t) = (-1)^{n\binom{m-1}{2}}(n+1)b^{mn^2+n}$$

$$D_2^{(m)}(mn+1,a,b,t) = (-1)^{n\binom{m-1}{2}}(t+(n+1)a)^2 b^{mn^2+3n}$$

$$D_2^{(m)}(mn-2,a,b,t) = -(-1)^{n\binom{m-1}{2}}(t+na)^2 b^{mn^2-3n}$$

$$D_2^{(m)}(mn-1,a,b,t) = (-1)^{n\binom{m-1}{2}}nb^{mn^2-n}$$

(7.10)

and $D_2^{(m)}(n,a,b,t) = 0$ *else*.

**Conjecture 7.7**

*Let*

$$G_m(z,a,b) = \frac{1}{\sqrt{(1-az)^2 - 4bz^m}} = \sum_{n\geq 0} g(n,m,a,b)z^n. \qquad (7.11)$$

*Then*

$$dd_0^{(m)}(mn,a,b) = (-1)^{n\binom{m-1}{2}} 2^{mn-1} b^{n(mn-1)},$$

$$dd_0^{(m)}(mn+1,a,b) = (-1)^{n\binom{m-1}{2}} 2^{mn} b^{n(mn+1)}$$

(7.12)

*and*

$dd_0^{(m)}(n,a,b) = 0$ *else.*

*For $m \geq 3$ we get*

$$dd_1^{(m)}(mn,a,b) = (-1)^{n\binom{m}{2}} 2^{mn} b^{mn^2},$$

$$dd_1^{(m)}(mn+1,a,b) = (-1)^{n\binom{m}{2}}(2n+1)2^{mn} ab^{mn^2+2n}$$

$$dd_1^{(m)}(mn-1,a,b) = -(-1)^{n\binom{m}{2}}(2n-1)2^{mn-2} ab^{mn^2-2n}$$

(7.13)

and $dd_1^{(m)}(n,a,b) = 0$ *else.*

*For $m \geq 4$ we get*

$$dd_2^{(m)}(mn,a,b) = (-1)^{n\binom{m-1}{2}}(2n+1)2^{mn}b^{n(mn+1)},$$

$$dd_2^{(m)}(mn+1,a,b) = (-1)^{n\binom{m-1}{2}}(2n+1)^2 2^{mn} a^2 b^{n(mn+3)},$$

$$dd_2^{(m)}(mn-2,a,b) = -(-1)^{n\binom{m-1}{2}}(2n-1)^2 2^{mn-3} a^2 b^{n(mn-3)},$$

$$dd_2^{(m)}(mn-1,a,b) = (-1)^{n\binom{m-1}{2}}(2n-1)2^{mn-1} b^{n(mn-1)}.$$

(7.14)

and $dd_2^{(m)}(n,a,b) = 0$ *else.*

**Remark 7.8**

Closely related is the following conjecture:

Let $m \geq 2$ and

$$\tilde{p}_n(m,z) := \det \begin{pmatrix} a(0) & a(1) & \cdots & a(n-1) & 1 \\ a(1) & a(2) & \cdots & a(n) & z \\ a(2) & a(3) & \cdots & a(n+1) & z^2 \\ \vdots & & & & \vdots \\ a(n) & a(n+1) & \cdots & a(2n-1) & z^n \end{pmatrix}$$

with $a(n) = G(n,m,a,b)$.

Then

$$\tilde{p}_{mn}(m,z) = (-1)^{n\binom{m-1}{2}} 2^{mn-1} b^{n(mn-1)} L_{2n}\left(\sqrt{z^{m-2}}(z-a), -b\right),$$

$$\tilde{p}_{mn+1}(m,z) = (-1)^{n\binom{m-1}{2}} 2^{mn} b^{n(mn+1)} z^{1-\frac{m}{2}} L_{2n+1}\left(\sqrt{z^{m-2}}(z-a), -b\right),$$ (7.15)

$$\tilde{p}_{mn+m-1}(m,z) = (-1)^{\binom{m-1}{2}} 2^{m-2} b^{(m-2)(2n+1)} z^{1-\frac{m}{2}} \tilde{p}_{mn+1}(m,z).$$

I want to conclude this paper with a further conjecture related to these results. We have

$$\frac{1}{\alpha(m,z) - \beta(m,z)} = \frac{1}{\alpha(m,z)} \left( \frac{1}{1 - \frac{\beta(m,z)}{\alpha(m,z)}} \right) = \frac{1}{\alpha(m,z)} \left( \frac{1}{1 - \frac{\alpha(m,z)\beta(m,z)}{\alpha(m,z)^2}} \right)$$

$$= \frac{1}{\alpha(m,z)} \left( \frac{1}{1 - \frac{bz^m}{\alpha(m,z)^2}} \right) = f_m(z) \left( \frac{1}{1 - bz^m f_m(z)^2} \right) = \sum_{n \geq 0} (bz^m)^n f_m(z)^{2n+1}.$$

Now consider the partial sums

$$H_m(k,z,a,b) = f_m(z,a,b)\frac{1-(bz^m f_m^2(z,a,b))^k}{1-bz^m f_m^2(z,a,b)} = G_m(z,a,b)\left(1-(bz^m f_m^2(z,a,b))^k\right).$$

Their coordinates coincide with those of $G_m(z,a,b)$ for $n < mk$.

Let $D^{(m,k)}(n,a,b)$ be the Hankel determinant of the coefficient sequence of $H_m(k,z,a,b)$.

Computer computations suggest that this sequence of Hankel determinants has a recurrent pattern modulo $mk$ which depends on the parity of $k$. More precisely we state the following

**Conjecture 7.10**

*Let $k, \ell \in \mathbb{N}$ and $k > 1$.*

*a1) For $mk\ell < n \leq mk(\ell+1)$ and $n \equiv jm \bmod(km)$*

$$D^{(m,2k)}(n,a,b) = (-1)^{n\binom{m-1}{2}}(\ell+2)^{jm-1}(\ell+1)^{m(k-j)}b^{n(nm-1)},$$

*a2) for $mk\ell \leq n < mk(\ell+1)$ and $n \equiv (jm+1)\bmod(km)$*

$$D^{(m,2k)}(n,a,b) = (-1)^{n\binom{m-1}{2}}(\ell+2)^{jm}(\ell+1)^{m(k-j)-1}b^{n(nm+1)}.$$

*b1) Let $m(2k+1)\ell < n \leq m(2k+1)(\ell+1)$. For $n \equiv jm \bmod(2k+1)m$ and $j \leq k$*

$$D^{(m,2k+1)}(n,a,b) = (-1)^{n\binom{m-1}{2}}(2\ell+2)^{jm-1}(2\ell+1)^{m(k-j)+1}b^{n(nm-1)},$$

*for $j \geq k+1$*

$$D^{(m,2k+1)}(n,a,b) = (-1)^{n\binom{m-1}{2}}(2\ell+3)^{(j-1-k)m}(2\ell+2)^{(2k+1-j)m}b^{n(nm-1)}.$$

*b2) Let $m(2k+1)\ell \leq n < m(2k+1)(\ell+1)$. For $n \equiv (jm+1)\bmod(2k+1)m$ and $j \leq k$*

$$D^{(m,2k+1)}(n,a,b) = (-1)^{n\binom{m-1}{2}}(2\ell+2)^{jm}(2\ell+1)^{m(k-j)}b^{n(nm+1)},$$

*for $j > k$*

$$D^{(m,2k+1)}(n,a,b) = (-1)^{n\binom{m-1}{2}}(2\ell+3)^{(j-k-1)m+1}(2\ell+2)^{m(2k+1-j)-1}b^{n(nm+1)}.$$

*All other Hankel determinants vanish.*

As an example consider the sequence of aerated Catalan numbers $(1,0,1,0,2,0,5,0,14,0,\cdots)$ with generating function $f_2(z,0,1) = \dfrac{1-\sqrt{1-4z^2}}{2z^2}$.

Their Hankel determinants are $d_0^{(2)}(n,0,1) = 1$ for all $n \in \mathbb{N}$. The Hankel determinants corresponding to
$$G_2(z,0,1) = \frac{1}{\sqrt{1-4z^2}} = \sum_{n\geq 0} \binom{2n}{n} z^{2n}$$ are $dd_0^{(2)}(n,0,1) = 2^{n-1}$ for $n \geq 1$.

By (2.21) $H_2(k,z,0,1) = \sum_{n\geq 0}\left(\binom{2n}{n} - \binom{2n}{n-k}\right) z^{2n}$.

It is interesting to see how the sequence of Hankel determinants $D^{(2,2k+1)}(n,0,1)$ converges to the sequence of Hankel determinants $dd_0^{(2)}(n,0,1)$.

$$\left(D^{(2,3)}(n,0,1)\right)_{n\geq 0} = \left(1,1,2 \mid 2^2, 2^2, 2\cdot 3 \mid 3^2, 3^2, 3\cdot 4 \mid 4^2, 4^2, 4\cdot 5 \mid \cdots\right),$$

$$\left(D^{(2,5)}(n,0,1)\right)_{n\geq 0} = \left(1,1,2,2^2,2^3 \mid 2^4, 2^4, 2^3\cdot 3, 2^2\cdot 3^2, 2\cdot 3^3 \mid 3^4, 3^4, 3^3\cdot 4, 3^2\cdot 4^2, 3\cdot 4^3 \mid \cdots\right),$$

$$\left(D^{(2,7)}(n,0,1)\right)_{n\geq 0} = \left(1,1,2,2^2,2^3,2^4,2^5 \mid 2^6, 2^6, 2^5\cdot 3, 2^4\cdot 3^2, 2^3\cdot 3^3, 2^2\cdot 3^4, 2\cdot 3^5, \mid 3^6, 3^6, 3^5\cdot 4, \cdots\right), \cdots.$$

Let me mention that the numbers $s(n)$ and $t(n)$ defined in (2.5) are $s(n,k) = 0$ and $t(kn,k) = \dfrac{n+2}{n+1}$, $t(kn-1,k) = \dfrac{n}{n+1}$ and $t(n,k) = 1$ else.